\documentclass[a4paper,preprint,12pt]{elsarticle} 
\usepackage[english]{babel} 
\usepackage[latin1]{inputenc}
\usepackage{amssymb}  
\usepackage{amsmath}    
\usepackage{pb-diagram}     
\usepackage{epstopdf}                   
\usepackage{graphicx}        
 \usepackage{epsfig}    
\usepackage{color}
\usepackage{fancyhdr}     
 \usepackage{natbib}  
\usepackage{rotating}
\usepackage{pdflscape}
\usepackage[T1]{fontenc}        
\usepackage{url}                
\usepackage{mathtools}                 
\usepackage{vmargin}    
\usepackage{lineno}
\usepackage{setspace}
\usepackage{color,soul}
\usepackage[ruled,vlined]{algorithm2e}  
\usepackage[title]{appendix} 
\biboptions{square}   
\usepackage[colorlinks=true,linkcolor= blue,urlcolor=blue,citecolor=blue]{hyperref} 
\usepackage[ruled,vlined]{algorithm2e} 
  
\newcommand{\bu}{\bullet}     
\newcommand{\bo}[1]{\mathbf{#1}} 

\def\indic{\hbox{1\kern-.24em\hbox{I}}}      

\newcommand{\var}{\mathbb{V}}  
\newcommand{\esp}{\mathbb{E}}    
\newcommand{\trace}{\text{Tr}}

\newcommand{\x}{x}
\newcommand{\X}{X}

\newcommand{\R}{\mathbb{R}}
\newcommand{\proba}{\mathbb{P}} 

\newcommand{\MG}{g} 
\newcommand{\BG}{\mathbf{G}}

\newcommand{\M}{f}

\newcommand{\T}{T}    
\newcommand{\N}{\mathbb{N}}      
\newcommand{\NN}{n}

\newcommand{\norme}[1]{\left|\left| #1 \right|\right|_{2}}    
\newcommand{\norml}[1]{\left|\left| #1 \right|\right|_{1}}

\newcommand{\normh}[1]{\left|\left| #1 \right|\right|_\mathcal{H}}

{\bf}{\it} 
\newtheorem{defi}{Definition}{\bf}{\it}  

\newtheorem{theorem}{Theorem}{\bf}{\it}     
\newtheorem{lemma}{Lemma}{\bf}{\it}          
\newtheorem{rem}{Remark}{\bf}{\it} 
\newtheorem{corollary}{Corollary}{\bf}{\it} 
\sethlcolor{green}      

\catcode`\"=\active 
\catcode`\"=\active 
\def"{\og\ignorespaces}
\def"{{\fg}} 
 
\usepackage{etoolbox}  

\makeatletter
\def\ps@pprintTitle{%
  \let\@oddhead\@empty
  \let\@evenhead\@empty
  \def\@oddfoot{\reset@font\hfil\thepage\hfil}
  \let\@evenfoot\@oddfoot
}
\makeatother

\begin{document}

\begin{frontmatter}  

\title{Kernel-based measures of association between inputs and outputs using ANOVA}
\author[a,b]{Matieyendou Lamboni\footnote{Corresponding author: matieyendou.lamboni[at]gmail.com/univ-guyane.fr, 19/11/2023}}         
\address[a]{University of Guyane, Department DFRST, 97346 Cayenne, French Guiana, France} 
\address[b]{228-UMR Espace-Dev, University of Guyane, University of R\'eunion, IRD, University of Montpellier, France.}    
                                  
\begin{abstract} 
ANOVA decomposition of function with random input variables provides ANOVA functionals (AFs), which contain information about the contributions of the input variables on the output variable(s). By embedding AFs into an appropriate reproducing kernel Hilbert space regarding their distributions, we propose an efficient statistical test of independence between the input variables and output variable(s). The resulting test statistic leads to new dependent measures of association between inputs and outputs that allow for i) dealing with any distribution of AFs, including the Cauchy distribution, ii) accounting for the necessary or desirable moments of AFs and  the interactions among the input variables. In uncertainty quantification for mathematical models, a number of existing measures are special cases of this framework. We then provide unified and general global sensitivity indices and their consistent estimators, including  asymptotic distributions. For Gaussian-distributed AFs, we obtain Sobol' indices and dependent generalized sensitivity indices using quadratic  kernels. 
\end{abstract}         
                                     
         
						  
\end{frontmatter}       
  
\section{Introduction} 
In statistical modeling and data analysis, symmetric measures of dependence for random vectors such as the Pearson correlation (\cite{pearson1901});  the Spearman correlation; the Kendall correlation; the canonical correlation coefficient (\cite{hotelling36,Kojadinovic09});  the RV coefficient for two random vectors (\cite{escoufier73}); the maximum mean discrepancy (\cite{borgwardt06,gretton07,gretton12}), including the energy distance (\cite{rizzo16}); the Hilbert-Schmidt independence criterion (\cite{gretton05,gretton05hsic,sejdinovic13}), including the  distance correlation (\cite{feuerverger93,szekely07}) and a recent correlation coefficient (\cite{chatterjee20}) have been used by different communities for testing whether two random vectors are independent or not, and then for measuring the dependence between such random vectors. Drawbacks, advantages and links between such dependent measures can be found in \cite{renyi59,sriperumbudur10,sejdinovic13,josse16} for instance.\\ 
    
The above dependent measures do not require any function linking both random vectors. In presence of a mathematical model or a black box function of the form $Y =\M(\bo{\X})$, where $\bo{\X} := (\X_1, \ldots, \X_d)$ are the model inputs with known probability distributions and $Y$ is the model output, such measures can still be used for testing the independence between the model output $Y=\M(\bo{\X})$ and some inputs such as $\X_1$. Such a problem occurs in computer experiments where a natural or human induced phenomenon is represented by complex computer code with  numerous input variables (\cite{derocquigny08}). For dimension reduction, it is relevant to identify unessential variables by means of different criteria such as statistical tests of independence or the variance-based sensitivity indices (SIs) (\cite{sobol93,saltelli00,lamboni11,gamboa14,lamboni18}), which were developed in the framework of the statistical theory called ANOVA (\cite{hoeffding48a,efron81,antoniadis84,sobol93}). \\ 
       
Formally, the above dependent measures between $\M(\bo{\X})$ and $\X_1$ rely on the test statistic of the form $T(\M(\bo{\X}), \X_1)$ or $T(\M(\bo{\X}), h(\X_1))$ with $T,\, h$ some given functions. When a statistical test reveals that the random variables $\M(\bo{\X})$ and $\X_1$ are dependent, the associated test statistic usually serves as a reasonable measure of the association between both variables. While this crude dependent measure is relevant for screening input variables in some cases, 
it may be not satisfying for ranking input variables because it involves only the inputs of interest and the output or the output and its conditional expectation given the inputs of interest (\cite{daveiga15,xiao18b,plischke20,barr22}). Moreover, it i) may result in dealing with transformations of the model output rather than using $Y$ (\cite{owen14}); ii) may include undesirable effects or complicate the problem (\cite{owen14,song12}); iii) ignores the asymmetric role of the model inputs and output.  For instance, the MMD-based indices result in using $\phi(\bo{Y})$ as outputs with $\phi$ some given functions a.k.a. feature maps (see \cite{daveiga21}). In the same sense, the HSIC-based indices proposed in \cite{daveiga21} result in working with transformations of the original inputs and output, as the HSIC requires defining one kernel on the inputs and another one on the outputs (\cite{smola07}).\\
                    
It is known that ANOVA remains an important step in exploratory and confirmatory analysis of model (\cite{gelman05}), and it relies on direct use of the model output(s) and inputs. The well-known ANOVA-like decomposition of $\M(\bo{\X})$ with independent variables $\bo{\X}$ can be written as follows (\cite{hoeffding48a,efron81,antoniadis84,sobol93,lamboni13,lamboni16}): 
\begin{equation} \label{eq:indeq} 
\M(\bo{\X}) =\sum_{\substack{v \subseteq \{1, \ldots, d\}}} \M_v\left(\bo{\X}_v\right) = \esp\left[\M(\bo{\X}) \right] +\M_1^{tot}(\bo{\X}) + \M_{\sim 1} (\X_2, \ldots, \X_d) \, ,   
\end{equation}                      
where $\bo{\X}_v$ is a subset of the model inputs whose subscripts belong to $v$; and $\M_1^{tot}(\bo{\X}) := \sum_{\substack{v \subseteq \{1, \ldots, d\}\\ v \cap \{1\} \neq \emptyset}} \M_v\left(\bo{\X}_v\right)$ denotes the total ANOVA functional (AF) of $\X_1$. The FANOVA decomposition of functions with non-independent variables is similar thanks to dependency models of dependent variables (see Section \ref{sec:SF} and \cite{lamboni21,lamboni21ar,lamboni22b}). ANOVA functionals are random variables, which contain the primary information about the contributions of the input variables on the output variable. In view of Equation (\ref{eq:indeq}),  $\M(\bo{\X})$ does not dependent on $\X_1$ whenever the total AF follows the Dirac probability measure $\delta_{\bo{0}}$, that is, 
$ 
\M_1^{tot}(\bo{\X}) = 0     
$ almost surely (a.s.) (see \cite{lamboni20} for independent variables and Lemma \ref{lem:cripr} in general). \\ 
                  
The aim of this paper is to propose a new dependent measure between the model inputs and output(s) that is more comprehensive and much flexible to account for the necessary or desirable statistical properties of AFs (e.g., higher-order moments, the Cauchy and heavy tailed distributions, asymmetric distributions) and the interactions among the input variables. The proposed dependent measure relies on the properties of AFs such as the independence criterion $\M_1^{tot}(\bo{\X}) =0 \; a.s.$ and the reproducing kernel Hilbert space (RKHS) (\cite{aronszajn50,Scholkopf02,berlinet04}). The RKHS theory offers a flexible framework for precise statistical inference of probability distributions (\cite{gretton07,smola07,sriperumbudur10}), and it is crucial for conducting independence tests that account for the necessary or sufficient or desirable statistical properties of distributions (\cite{gretton05,gretton07,smola07,fukumizu08,fukumizu09,song12}). The proposed dependent measure leads to provide the first-order and total kernel-based sensitivity indices, and its empirical expression is used for deriving a test statistic for testing the independence between the model inputs and output(s).\\       
                                 
This paper is organized as follows: Section \ref{sec:SF} deals with AFs and their statistical properties in terms of their ability to assess the effects of the input variables on the model output(s).  Such properties enable the formulation of the initial null hypothesis for the independence test between input variables and the model output(s) and equivalent null hypothesis such as the variance   $\var\left[\M_1^{tot}(\bo{\X})\right]=0$. Most of variance-based SIs rely on that null hypothesis for performing dimension reduction. Indeed, $\var\left[\M_1^{tot}(\bo{\X})\right]$ represents the non-normalized total SI of $\X_1$.\\    
Since the associated alternative hypothesis is sufficient for AFs that are Gaussian-distributed and the Cauchy distribution does not have the variance, Section \ref{sec:equivemb} is devoted to develop equivalent null hypothesis for the independence test by embedding AFs into an appropriate RKHS regarding their distributions.  Moreover, as a given kernel can help for working with specific moments (\cite{song12}), we also provide the set of distribution-free kernels that guarantee the equivalent null hypothesis for the independence test between the model inputs and output(s). Section \ref{sec:teststat} presents a statistical test that leads to an interesting measure of association between inputs and output(s). In Section \ref{sec:empkbsi}, we formally introduce our empirical dependent measure and the associated test statistic. We also provide and study kernel-based SIs.  We provide analytical and numerical results in Section \ref{sec:test} and conclude this work in Section~\ref{sec:con}.      
                                                             
\section*{General notations}    
For an inetegr $d>1$, $d$ input variables $\bo{\X} :=(X_1, \, \ldots,\, X_d)$ and $u \subseteq \{1,\,\ldots,\, d\}$, we use $\bo{\X}_u :=(\X_j, \forall\, j \in u)$, $\bo{\X}_{\sim u} :=(\X_j, \forall\, j\in \{1,\,\ldots,\, d\}\setminus u)$ and $|u|$ for the number of elements in $u$. Thus, we have the partition $\bo{\X} =(\bo{\X}_{u},\, \bo{\X}_{\sim u})$. We also use $\bo{R} \stackrel{d}{=} \bo{\X}$ to say that $\bo{R}$ and $\bo{\X}$ have the same cumulative distribution function (CDF).
             
\section{Theoretical properties of ANOVA functionals of (in)-dependent inputs} \label{sec:SF}
In this section, we provide the properties of ANOVA functionals for functions with non-independent variables. Such properties serve as the initial null hypothesis for testing the independence between the random input variables and output variables. For the sequel of generality, consider a vector-valued function  $\M : \R^d \to \R^\NN$, which includes a $d$-dimensional random vector $\bo{\X}$ as inputs and provides $\NN$ outputs given by $\M(\bo{\X})$. We are interested in measuring the association between the following two random vectors: $\M(\bo{\X})$ and $\bo{\X}_u$ with $u \subseteq \{1, \ldots, d\}$ and $u\neq \emptyset$, keeping in mind the asymmetric role between both random vectors.  \\    
    
For any distribution of $\bo{\X}$, we are able to model $\bo{\X}$ as follows (\cite{skorohod76,lamboni21,lamboni22b,lamboni21ar}):             
\begin{equation}   \label{eq:dep} 
\bo{\X}_{\sim u} \stackrel{d}{=} r_u\left(\bo{\X}_{u},  \bo{Z} \right)\, ; 
\end{equation} 
where $r_u$ is a function;  $\bo{Z}$ is a random vector of $d-|u|$ independent variables, and $\bo{\X}_{u}$ is independent of $\bo{Z}$. Note that when $\bo{\X}$ is consisted of independent variables, the dependency function in (\ref{eq:dep}) comes down to $\bo{\X}_{\sim u} \stackrel{d}{:=} r_u\left(\bo{Z} \right)$. 
Composing $\M$ by the function $r_u$ in (\ref{eq:dep}) yields    
\begin{equation} \label{eq:equi}  
\MG(\bo{\X}_{u}, \bo{Z}) \stackrel{d}{:=} \M(\bo{\X}_{u}, r_u(\bo{\X}_{u}, \, \bo{Z})) \, . 
\end{equation}  
 In view of Equation (\ref{eq:equi}), the function $\MG$ includes two independent random vectors (i.e., $\bo{\X}_{u}$ and $\bo{Z}$) as inputs, and it provides $\NN$ outputs sharing the distribution of $\M(\bo{\X})$. For independent variables $\bo{\X}$, we can see that  
$ \MG(\bo{\X}_{u}, \bo{Z}) = \M(\bo{\X}_{u}, r_u(\bo{Z})) \stackrel{d}{=} \M(\bo{\X}_{u}, \bo{\X}_{\sim u})$. For the sequel of generality, we are going to use $\MG(\bo{\X}_{u}, \bo{Z})$ in what follows,  which comes down to $\M(\bo{\X}_{u}, \bo{\X}_{\sim u})$ when the inputs are independent.  As a matter of fact, we
can always claim that \\  
  
\noindent 
(A1) our function of interest $g$ includes two independent random vectors: $\bo{\X}_{u}$ and $\bo{Z}$. \\
         
Note that (A1) is not theoretically or practically restrictive because it is always satisfied thanks to dependency models (see (\ref{eq:dep})).  
          
\subsection{ANOVA functionals}  
Under (A1), the Hoeffding decomposition of $\MG(\bo{\X}_{u}, \bo{Z})$ is given by (\cite{hoeffding48a,efron81,antoniadis84,sobol93})   
$$ 
\MG(\bo{\X}_{u}, \bo{Z}) =  \esp_{\bo{\X}_{u},\bo{Z}}\left[\MG(\bo{\X}_{u}, \bo{Z})\right] + 
\MG_u(\bo{\X}_{u}) + \MG_{\sim u}(\bo{Z}) + \MG_{u,\sim u}(\bo{\X}_{u}, \bo{Z})\, ,    
$$    
where  $\esp_{\bo{\X}_{u},\bo{Z}}\left[\MG(\bo{\X}_{u}, \bo{Z})\right] =: \esp\left[\MG(\bo{\X}_{u}, \bo{Z})\right]$ is the expectation taking w.r.t. $\bo{\X}_{u},\bo{Z}$;    
$$
\MG_u(\bo{\X}_{u}) :=  \esp_{\bo{Z}}\left[ \MG(\bo{\X}_{u}, \bo{Z})\right] - \esp\left[\MG(\bo{\X}_{u}, \bo{Z})\right];  \,
\quad 
\MG_{\sim u} (\bo{Z}) :=  \esp_{\bo{\X}_u}\left[\MG(\bo{\X}_{u}, \bo{Z})\right] - \esp\left[\MG(\bo{\X}_{u}, \bo{Z})\right] \, ,   
$$ 
$$ 
\MG_{u,\sim u}(\bo{\X}_{u}, \bo{Z}) := \MG(\bo{\X}_{u}, \bo{Z}) - \MG_u(\bo{\X}_{u}) - \MG_{\sim u}(\bo{Z}) + \esp\left[\MG(\bo{\X}_{u}, \bo{Z})\right] \, .   
$$  
 
The first-order and total AFs of $\bo{\X}_u$ are defined as follows: (\cite{lamboni18,lamboni19,lamboni21,lamboni21ar,lamboni22,lamboni22b})  
$$
\MG^{fo}_u(\bo{\X}_{u}) := \esp_{\bo{Z}}\left[ \MG(\bo{\X}_{u}, \bo{Z})\right] - \esp\left[\MG(\bo{\X}_{u}, \bo{Z})\right] = \MG_u(\bo{\X}_{u})  \, ,     
$$        
$$
\MG^{tot}_u (\bo{\X}_{u}, \bo{Z}) := \MG(\bo{\X}_{u}, \bo{Z}) - \esp_{\bo{\X}_{u}} \left[\MG(\bo{\X}_{u}, \bo{Z}) \right] \, .    
$$   
It is worth noting that both AFs are zero-mean and $\NN$-dimensional random vectors, which are directly based on the model outputs. We also have the following relationship:  
\begin{equation} \label{eq:linktotfo}
\MG^{fo}_u (\bo{\X}_{u}) =   \esp_{\bo{Z}}\left[\MG^{tot}_u (\bo{\X}_{u}, \bo{Z}) \right] \, .      
\end{equation} 
Recall that AFs contain the primary information about the contribution of the random vector $\bo{\X}_{u}$ to the model outputs given by $\MG(\bo{\X}_{u}, \bo{Z}) \stackrel{d}{=} \M(\bo{\X})$. The main interesting property of AFs is derived in  Lemma \ref{lem:cripr}.         
\begin{lemma}  \label{lem:cripr}    
Under (A1), $\M(\bo{\X}) \stackrel{d}{=} \MG(\bo{\X}_{u}, \bo{Z})$ is independent of $\bo{\X}_u$ iff   
$$
 \MG^{tot}_u (\bo{\X}_{u}, \bo{Z}) = \bo{0} \; \;  a.s. \, .            
$$      
\end{lemma}       
\begin{preuve}     
See Appendix \ref{app:lem:cripr}. 
\begin{flushright} 
$\Box$ 
\end{flushright} 
\end{preuve}    

From Lemma  \ref{lem:cripr}, it is clear that the total AF $\MG^{tot}_u(\bo{\X}_u, \bo{Z})$ is sufficient for fully characterizing the independence between $\M(\bo{\X})$ and $\bo{\X}_u$. We then formulate the initial test hypotheses of independence as follows:   
$$
H_0: \; \MG^{tot}_u(\bo{\X}_u, \bo{Z}) \stackrel{}{=} \bo{0} \; a.s. \, 
\qquad \mbox{Vs.} \qquad   
H_1: \; \MG^{tot}_u(\bo{\X}_u, \bo{Z}) \stackrel{}{\neq} \bo{0} \; a.s. \, .   
$$   
While the null hypothesis is equivalent to   
$  
 \var\left[\MG^{tot}_u(\bo{\X}_u, \bo{Z}) \right] =\mathsf{O}
$,  
 with $\mathsf{O}$ a null matrix, its alternative hypothesis given by $ \var\left[\MG^{tot}_u(\bo{\X}_u, \bo{Z}) \right] \neq \mathsf{O}$ may be not satisfying because it does not account for higher-order moments. Moreover, this null hypothesis implicitly requires the existence of the variance-covariance of AFs, which is not the case for the Cauchy distribution for instance. Therefore, we need a dependent measure that accounts for sufficient or desirable higher-order moments for a given distribution of $\MG^{tot}_u(\bo{\X}_u, \bo{Z})$. Since the total AF is a random vector, we may use the probability metrics based on the difference between distribution functions or the Wasserstein metric or the kernel methods for making such comparisons. Significant discussions about such metrics and kernel methods can be found in \cite{sriperumbudur10}. To develop our dependent measure, we are going to use the kernel methods that are much flexible for including specific moments of AFs. 
              
\subsection{Embedding ANOVA functionals into a RKHS}    
To build a statistical test that i) leads to interesting dependent measures, ii) is much flexible for explicitly including specific moments, iii) is able to distinguish two different distributions, AFs are going to be embedded into a RKHS or feature spaces according to their distributions (\cite{berlinet04,gretton07,smola07}).  
  
\begin{defi} (Aronszajn, 1950)          
Let $\mathcal{X}$ be an arbitrary space, $\mathcal{H}$ be an Hilbert space  endowed with the inner product $\left< \cdot,\, \cdot \right>$. The functions\\  
  
(i) $\phi : \mathcal{X} \to \mathcal{H}$ is called a feature map;\\      
  
(ii) $k : \mathcal{X} \times \mathcal{X} \to \R$ given by $k(r, r') = \left< \phi(r),\, \phi(r') \right> =: \left<k(\cdot, r),\, k(\cdot, r') \right>$ is called a valid kernel. 
\end{defi}     
A kernel $k \in \mathcal{K}$ is said centered at $r_0$ when
$k(\cdot, r_0) =0$.  Given a kernel $k(r, r')$, we can construct a new kernel that is centered at $r_0$ as follows (\cite{sejdinovic13}): 
$$
k_c\left(r, r' \right) := k\left(r, r' \right) + k\left(r_0, r_0 \right) - k\left(r, r_0 \right) - k\left(r_0, r' \right)\, .  
$$   
  
For a $\NN$-dimensional random vector $\bo{G}$ having $F$ as CDF (i.e., $\bo{G} \sim F$) such as AFs, the transformation $k(\cdot, \bo{G})$ aims at embedding this random vector into a RKHS induced by $k(\bo{G}, \bo{G}')$ with $\bo{G}'$ an i.i.d. copy of $\bo{G}$.    
Linear statistics in the new RKHS such as the mean element account for all the moments or desirable moments of $\bo{G}$ depending on the kernel  $k$ (\cite{gretton07,smola07,gretton12}). For embedding AFs into a RKHS and working with the mean element (see Definition \ref{def:meanel}), we use $\mathcal{P}_B$ for the set of Borel probability distributions and define the set of AFs distributions as follows:      
$$  
\mathcal{F} := \left\{F \in \mathcal{P}_B : \, \MG^{tot}_u (\bo{\X}_{u}, \bo{Z}) \sim F \quad  \mbox{or}\quad  \MG^{fo}_u (\bo{\X}_{u}) \sim F,\; \;  \forall \, u \subseteq \{1, \ldots, d\}  \right\}\, . 
$$      
 
The class of distributions $\mathcal{F}$ is adequate for manipulating AFs. For a valid and measurable kernel, we assume that \\ 

(A2): $\esp\left[\sqrt{k(\bo{G}, \bo{G})} \right] <~\infty$ for all $\bo{G}\sim F \in \mathcal{F}$.  
         
\begin{defi} (\cite{gretton07,gretton12}).  \label{def:meanel} 
 Consider a kernel $k$ and $\bo{G}\sim F \in \mathcal{F}$. \\ 
The map $\bo{G} \mapsto \esp_F\left[k(\cdot, \bo{G}) \right]$ is called a mean feature map; and $\mu_{F}(\bo{G}) := \esp_F\left[k(\cdot, \bo{G}) \right]$ is called the mean element. 
\end{defi} 
         
Generally, the feature map is used to embed $\bo{G}$ into a higher dimensional random vectors so as to include all type of information about the data. Characteristic kernels aim to accomplish such tasks (\cite{fukumizu04,gretton07,fukumizu08,fukumizu09,fukumizu09b}).   
   
\begin{defi} (\cite{gretton07,fukumizu08}) \label{def:chker} 
Conisder two random vectors $\bo{G}_1 \sim F, \; \bo{G}_2 \sim H$. \\       
\noindent 
A measurable kernel $k$ is said characteristic if the mean feature map is one to one:
$$   
\mu_{F}(\bo{G}_1) = \mu_{H}(\bo{G}_2) \Longleftrightarrow  F= H \quad (i.e., \bo{G}_1 \stackrel{d}{=} \bo{G}_2) \, .          
$$  
\end{defi} 
Taking the distance between $\mu_{F}(\bo{G}_1)$ and $\mu_{H}(\bo{G}_2)$ leads to the maximum mean discrepancy (MMD), that is,  
$
MMD^2(\bo{G}_1, \bo{G}_2) := \normh{\mu_{F}(\bo{G}_1) - \mu_{H}(\bo{G}_2)}^2
$ (\cite{borgwardt06,gretton07,gretton12}).  Note that the centered kernel $k_c$ associated with a valid kernel $k$ is a characteristic kernel if and only if $k$ is a characteristic one (\cite{sejdinovic13}). 
\\ 
    
Thus, the mean element associated with a characteristic kernel uniquely determines a probability distribution. This interesting property gives us the ability to use the mean element for fully characterizing the distribution of AFs. It is to be noted that characteristic kernels for specific class of distributions such as the class of Gaussian distributions can be defined as well. Thus, characteristic kernels on $\mathcal{F}$ are sufficient for distinguishing different AFs.     
 For instance, while a test statistic that can distinguish the first two moments is sufficient for a class of Gaussian distributions, the mean element of the form $\esp\left[e^{\bo{G}^\T t}\right]$ with $t \in \R^\NN$, which generalizes the notion of moment-generating function in probability, allows for distinguishing all the moments of a probability distribution. 
     
\section{Kernels for equivalent null hypothesis} \label{sec:equivemb}
This section provides a set of kernels that ensure the equivalent criterion of independence between the input variables $\bo{\X}_u$ and the model outputs. 
Recall that the null hypothesis $\var\left[\MG^{tot}_u(\bo{\X}_u, \bo{Z})\right] =\mathsf{O}$ is an equivalent criterion of independence, and it is also obtained using the quadratic kernel, that is, $k_2\left(\bo{r}, \bo{r}'\right) =\left< \bo{r}, \, \bo{r}'\right>^2_{\R^\NN}$. Indeed, if we use $\bo{\X}_u', \bo{Z}'$ for i.i.d. copies of $\bo{\X}_u, \bo{Z}$, we can check that  
$$ 
 \esp \left[k_2\left( \MG^{tot}_u(\bo{\X}_u, \bo{Z}), \, \MG^{tot}_u(\bo{\X}_u', \bo{Z}') \right) \right] = 0  \Longleftrightarrow  \var\left[\MG^{tot}_u(\bo{\X}_u, \bo{Z})\right] =\mathsf{O}\, .  
$$  
    
It is worth noting that kernels of the form $k_{2p}\left(\bo{r}, \bo{r}'\right) =\left< \bo{r}, \, \bo{r}'\right>^{2p}_{\R^\NN}$ for integer $p >0$ lead to an equivalent null hypothesis of independence, although such kernels are not characteristic in general. Since some kernels do not ensure the independence criterion, let us start with the following definitions thanks to Lemma \ref{lem:cripr}. Namely, we use $\mathcal{K}$ for the set of valid and measurable kernels; $H$ for the CDF of the Dirac probability measure $\delta_{\{\bo{0}\}}$.    
       
\begin{defi}   \label{def:equitest}        
A kernel $k \in \mathcal{K}$ is said to be an equivalent kernel for the independence test between $\bo{\X}_u$ and $\MG(\bo{\X}_{u}, \bo{Z})$  whenever 
\begin{equation} \label{eq:indcriker01}     
 \esp\left[k\left(\MG^{tot}_u (\bo{\X}_{u}, \bo{Z}), \,  \MG^{tot}_u (\bo{\X}_{u}', \bo{Z}') \right)\right] -k(\bo{0}, \bo{0})  = 0  \,  
\Longrightarrow    \,  
\MG^{tot}_u (\bo{\X}_{u}, \bo{Z}) = \bo{0} \; a.s. \; \, .  
\end{equation}
\end{defi}      
   
The equivalence used in Definition \ref{def:equitest} is guaranteed by Lemma \ref{lem:cripr}. For a centered kernel at $\bo{0}$ (i.e., $\bar{k}$), the left term of Equation (\ref{eq:indcriker01}) becomes 
$
\esp\left[\bar{k}\left(\MG^{tot}_u (\bo{\X}_{u}, \bo{Z}), \,  \MG^{tot}_u (\bo{\X}_{u}', \bo{Z}') \right)\right]=~0
$.  \\

To construct the set of equivalent kernels for the independence test, consider the following set of kernels:    
\begin{equation} \label{eq:setker} 
\mathcal{K}_E :=\left\{
k \in \mathcal{K} : \,  
\begin{array}{c}  
  \, \esp_{\bo{G} \sim F, \bo{G}'\sim F}\left[\bar{k}(\bo{G}, \bo{G}')\right] =0 \Rightarrow F= H ,\; \forall \, F \in \mathcal{F} \; \\
\mbox{or} \\
\esp_{(\bo{G},\, \bo{G}') \sim \nu}\left[k(\bo{G}, \bo{G}')\right] =0 \Rightarrow \nu=  0, 
 \; \forall\, \nu = F\otimes F- H\otimes H \,   
\end{array}   
\right\}     \, .             
\end{equation}  
        
We can check that $\mathcal{K}_E$ contains quadratic kernels that are centered at $\bo{0}$, and we are going to see that it contains some well-known characteristic kernels. Lemma \ref{lem:extkerind} gives interesting properties of $\mathcal{K}_E$ regarding the null hypothesis. 
 
\begin{lemma} \label{lem:extkerind} Let $k \in \mathcal{K}_E$ and assume that (A1) and (A2) hold.  Then,\\                        
  
$k$ is an equivalent kernel for the independence test between $\bo{\X}_u$ and $\MG(\bo{\X}_{u}, \bo{Z})$. 
\end{lemma}       
\begin{preuve} 
See Appendix \ref{app:lem:extkerind}. 
\begin{flushright}          
$\Box$          
\end{flushright}
\end{preuve}  
  
Lemma \ref{lem:extkerind} shows that the set $\mathcal{K}_E$ given by (\ref{eq:setker}) contains equivalent kernels for the independence criterion whatever are the distributions of the model outputs and total AFs. Despite $\mathcal{K}_E$ is not exhaustive, Lemma \ref{lem:kerequiv} shows its richness. 
To that end, consider the famous radial-based characteristic kernels on $\R^\NN \times \R^\NN$ given by 
$$ 
k(\bo{r}, \bo{r}') := \int e^{-i (\bo{r}-\bo{r}')^\T \bo{w}} \, d\Lambda(\bo{w}) \, , 
$$        
where $\Lambda$ is a positive and bounded Borel measure with the support $Supp(\Lambda) =\R^\NN$.  

\begin{lemma} \label{lem:kerequiv} Let $k,\, \bar{k} \in \mathcal{K}$ and assume (A1)-(A2) hold. \\ 

$\quad$ (i) If $k$ is a characteristic kernel, then the centered kernel $\bar{k} \in \mathcal{K}_E$. \\

 $\quad$ (ii) If $k$ is a radial-based characteristic kernel, then $k \in \mathcal{K}_E$.
\end{lemma}  
\begin{preuve} 
See Appendix \ref{app:lem:kerequiv}.   
\begin{flushright}     
$\Box$             
\end{flushright} 
\end{preuve}       
 
Kernels $\bar{k}$ of Point (i) in Lemma \ref{lem:kerequiv} lead to a comparison between the distributions of the total AF and the Dirac measure using the maximum mean discrepancy (see Section \ref{sec;linkmmd}). Thus, Point (ii) offers other possibilities that allow for working with non-centered and characteristic kernels such as Gaussian kernels. \\ 
        
\noindent      
\textbf{Example of equivalent kernels for the independence test} 
\begin{itemize} 
\item A class of distance-induced characteristic kernels that are already centered at $\bo{0}$ contains the following kernels (\cite{sejdinovic13}):    
$$
k_d^\alpha(\bo{r}, \bo{r}') := \frac{1}{2}\left( \norme{\bo{r}}^\alpha + \norme{\bo{r}'}^\alpha - \norme{\bo{r} -\bo{r}'}^\alpha \right), \quad \forall\, \alpha \in ]0,\, 2[ \, .    
$$        
  
\item Recall that kernels given by $k_{2p}(\bo{r}, \bo{r}') := \left< \bo{r}, \bo{r}' \right>^{2p}_{\R^\NN}$ for any integer $p \geq 1$ belong to $\mathcal{K}_E$. Moreover,  for any integer $2 \leq L < \infty$ and positive definite and diagonal matrix $\mathcal{D}_d$,  kernels of the form   
$
\bar{k}_{\mathcal{D}}(\bo{r}, \bo{r'}) := \sum_{q=1}^{L\geq 2} \left(\bo{r}^\T \mathcal{D}_q \bo{r}'\right)^{q}     
$           
are in $\mathcal{K}_E$. Each kernel $k_{q}\left(\bo{r}, \bo{r}' \right) := \left(\bo{r}^\T \mathcal{D}_q \bo{r}'\right)^q$ is sufficient for incorporating the $q^{th}$-order moments of the total AF and all the correlations among the components of AFs. In general, we are able to incorporate all the moments by taking $L \to +\infty$, which leads to the exponential kernel, that is, $k_{e} (\bo{r}, \bo{r'}) = e^{\alpha \left< \bo{r}, \, \bo{r}' \right>}$ with $\alpha >0$.  
     
\item Finally, radial-based characteristic kernels such as Gaussian, Laplacian and  the Cauchy kernels and their associated centered kernels belong to $\mathcal{K}_E$. Note that the Gaussian kernel is the normalized version of the exponential kernel. 
\end{itemize}     
 
\section{New independence test and dependent measure between inputs and outputs} \label{sec:teststat}
For testing the independence between $\MG\left(\bo{\X}_u, \bo{Z} \right)$ and $\bo{\X}_u$ and obtaining an interesting dependent measure when $\bo{\X}_u$ contribute to the outputs $\MG\left(\bo{\X}_u, \bo{Z}\right)$, we are going to use kernels that  guarantee the equivalent null hypothesis of independence between these random vectors such as the set of kernels $\mathcal{K}_E$.

\subsection{Test hypotheses and deviation mesaure from independence} \label{sec:dev}
For concise notations, we use $\BG^{tot}_u  \stackrel{}{:=} \MG^{tot}_u\left(\bo{\X}_u, \bo{Z} \right)$, $\BG^{tot '}_u \stackrel{}{:=} \MG^{tot}_u\left(\bo{\X}_u', \bo{Z}'\right)$, and we see that $\BG^{tot '}_u$ is an i.i.d. copy of $\BG^{tot}_u$. For a real $q>0$, $k\in \mathcal{K}_E$,  the generic test hypotheses are formally given by  
$$   
H_0' : \left|\esp \left[k \left(\BG^{tot}_u, \BG^{tot '}_u \right) \right] - k(\bo{0}, \bo{0}) \right|^q =0 
\quad \mbox{Vs}  \quad  
H_1' :  \left| \esp \left[k \left(\BG^{tot}_u, \BG^{tot '}_u \right) \right] - k(\bo{0}, \bo{0})  \right|^q \neq 0 \,  .    
$$         

If we use $F_{T_u} \in \mathcal{F}$  for the CDF of $\BG^{tot}_u$, we can measure the deviation from independence as follows: 
\begin{equation}  \label{eq:depm} 
\mathcal{D}_k^q(F_{T_u}) := \left| \esp_{F_{T_u}}\left[k\left(\BG^{tot}_u, \, \BG^{tot '}_u \right)\right] - k(\bo{0}, \bo{0}) \right|^q, \quad \forall\, q>0 \, .     
\end{equation}           
When $q=1$, $\mathcal{D}_k(F_{T_u})$ stands for $\mathcal{D}_k^1(F_{T_u})$. The discrepancy measure in (\ref{eq:depm} ) is still valid for any CDF $F \in \mathcal{F}$ such as the CDF $F_u$ of the first-order  AFs, that is, $\BG^{fo}_u  \stackrel{}{:=} \MG^{fo}_u\left(\bo{\X}_u\right) \sim F_u$.  A reasonable measure of the deviation from independence (or equivalently kernel) must be able to account for the fact that the first-order AFs bring partial information compared to the total ones (see Equation (\ref{eq:linktotfo})). This leads to the following definition.              
\begin{defi}  \label{def:scker} 
Consider AFs given by $\BG^{fo}_u \sim F_{u}$, $\BG^{tot}_u \sim F_{T_u}$ and a kernel $k \in  \mathcal{K}$.\\
  
 The kernel $k$ is said to be ANOVA compatible whenever   
$$ 
\mathcal{D}_k(F_{u}) \leq \mathcal{D}_k(F_{T_u}), \quad \forall \, u \subseteq \{1, \ldots, d\}   \, .
$$   
\end{defi}       

Using Jensen's inequality, we can check that the quadratic kernel is ANOVA compatible while the Hellinger kernel given by $k_H(r, r') := \sqrt{r r'}$ is clearly not. Combining the notion of ANOVA-compatible kernels and equivalent kernels for the independence test leads to the definition of importance-measure kernels (IMKs).  
\begin{defi}  \label{def:imk} 
A valid kernel $k$ is said to be an IMK whenever $k$ is ANOVA compatible and $k \in  \mathcal{K}_E$. We use  $\mathcal{K}_{IM}$ for the set of IMKs.   
\end{defi}    
          
While the quadratic kernel $\left< \bo{r}, \bo{r}' \right>^2_{\R^\NN}$ is part of $\mathcal{K}_{IM}$, it is not a characteristic kernel in general. Indeed, it is characteristic kernel for the class of Gaussian distributed AFs. Lemma \ref{lem:csak} provides some conditions for kernels to be IMKs. 
\begin{lemma} \label{lem:csak} Let $k, \,\bar{k} \in \mathcal{K}_{E}$ be kernels, $c \in \R$ and assume (A1)-(A2) hold.\\     
     
 $\quad$ (i) If $k(\bo{r}, \bo{0}) =c$ and $k(\bo{r}, \bo{r}')$ is convex in $\bo{r}$, then $\bar{k} \in \mathcal{K}_{IM}$. \\      
  
$\quad$ (ii) If $k(\bo{r}, \bo{r}')$ is convex in $\bo{r}$, then $k \in \mathcal{K}_{IM}$. \\ 

$\quad$ (iii) If $k(\bo{r}, \bo{r}')$ is concave in $\bo{r}$ and $k(\bo{0}, \bo{0}) > 0$, then $k \in \mathcal{K}_{IM}$.
\end{lemma}       
\begin{preuve}   
See Appendix \ref{app:lem:csak}. 
\begin{flushright}
$\Box$  
\end{flushright}
\end{preuve}    
 
From Lemma \ref{lem:csak}, convex and some concave kernels are IMKs. Of course, the assumptions of convexity or concavity are required on the support of the output distribution. For log-concave kernels of the form $\exp\left(-\alpha \psi(\bo{r}, \bo{r}') \right)$ with $\psi$ a convex function and $\alpha >0$, we are able to control such kernels through $\alpha$ in order to obtain concave kernels  on the support of the outputs (see Section \ref{app:conkerp}). \\  
             
\noindent       
\textbf{Examples of IMKs on $\R^d$.} 
\begin{itemize}
\item The quadratic kernel of the form $k_2(\bo{r}, \bo{r}') :=  \left(\bo{r}^\T\Sigma \bo{r}' \right)^2$ with $\Sigma$ a diagonal and positive definite matrix.
\item The absolute kernel or the $L_1$-based kernel of the form $k_a(\bo{r}, \bo{r}') :=  \norml{\bo{r}} \norml{\bo{r}'}$.  
\item The moment-generating kernel or the exponential kernel given by  $k_e(\bo{r}, \bo{r}') = \exp(\alpha \left<\bo{r}, \bo{r}' \right>)$ and its associated centered kernel at $\delta_0$ with $\alpha \in \R_+$. 
\item  The Laplacian kernel and the Gaussian kernel for some values of $\alpha >0$ (see Section \ref{app:conkerp}).   
\end{itemize}   

The IMKs provided in Lemma \ref{lem:csak} ensure interesting properties of the  discrepancy measure $\mathcal{D}_k^q(F)$ defined in Equation (\ref{eq:depm}). To provide such results in Theorem~\ref{theo:sckerg}, let us consider the input variables  $\bo{\X}_w$ with  $w \subseteq u$. We use $F_w$ (resp. $F_{T_w}$) for the CDF of the first-order (resp. total) AF of $\bo{\X}_w$, that is, $\BG^{fo}_w \sim F_{w}$ and $\BG^{tot}_w \sim F_{T_w}$.  
                       
\begin{theorem} \label{theo:sckerg} 
Let $k$ be an IMK given in Lemma \ref{lem:csak}; $w \subseteq u \subseteq \{1, \ldots, d\}$, $\BG^{fo}_u \sim F_{u}$ and $\BG^{tot}_u \sim F_{T_u}$ be AFs. Assume (A1)-(A2) hold. Then,  \\
     
\begin{equation}    
\mathcal{D}_k(F_{w}) \leq   \mathcal{D}_k(F_{u})  \, ;
\end{equation} 
\begin{equation}  
\mathcal{D}_k(F_{T_w}) \leq \mathcal{D}_k(F_{T_u})  \, .  
\end{equation} 
\end{theorem}  
\begin{preuve}  
See Appendix \ref{app:theo:sckerg}. 
\begin{flushright}        
$\Box$       
\end{flushright} 
\end{preuve}

It comes out from Theorem \ref{theo:sckerg} that the discrepancy measure increases with the cardinality of a subset of inputs. The fact that increasing the number
of components in a subset of inputs does not make the discrepancy measure smaller is commonly expected, as that property is satisfied in ANOVA and variance-based sensitivity analysis. Thus, Lemma \ref{lem:csak} provides kernels that guarantee interesting properties encountered in ANOVA. 
          
\subsection{Test statistic and dependent measure between inputs and outputs} \label{sec;impmea}
Based on the results from Theorem \ref{theo:sckerg}, it becomes clear that IMKs given in Lemma \ref{lem:csak} and Equation (\ref{eq:depm}) can lead to a coherent dependence measure of association between inputs and outputs according to R\'enyi' axioms (\cite{renyi59}). To define such dependent measures, we use $F_\bu \in \mathcal{F}$ for the CDF of the centered outputs $\bo{Y}:=\MG\left(\bo{\X}_u, \bo{Z} \right)- \esp[\MG\left(\bo{\X}_u, \bo{Z} \right)]$ and $\bo{Y}'$ for an i.i.d. copy of $\bo{Y}$.  
      
\begin{defi} \label{def:kbimnor}
For a real $q >0$, let $k \in \mathcal{K}_{IM}$ be an IMK and $F \in \mathcal{F}$ be the CDF of any AF.\\  
 The dependent measure of a random vector $\bo{G} \sim F$ is defined by 
\begin{equation}  \label{eq:depmesn} 
S^{k,q}_F := \frac{\mathcal{D}_k^q(F)}{\mathcal{D}_k^q(F_\bu)} = \left| \frac{\esp_{\bo{G} \sim F, \bo{G}' \sim F} \left[ k(\bo{G}, \bo{G}') \right] - k(\bo{0}, \bo{0})}{\esp_{\bo{Y} \sim F_\bu, \bo{Y}' \sim F_\bu} \left[ k(\bo{Y}, \bo{Y}') \right] - k(\bo{0}, \bo{0})} \right|^q  \, .   
\end{equation}      
\end{defi}            
               
For IMKs provided in Lemma \ref{lem:csak}, the right term of Equation (\ref{eq:depmesn})  can be written without the absolute symbol. Moreover, the dependent measures of the first-order and total AFs of $\bo{\X}_u$ having $F_u$ and $F_{T_u}$ as CDFs are given by, respectively  
$$ 
S^{k,q}_{F_{u}} := \frac{\mathcal{D}_k^q(F_{u})}{\mathcal{D}_k^q(F_\bu)};
\qquad  \qquad 
S^{k,q}_{F_{T_u}} := \frac{\mathcal{D}_k^q(F_{T_u})}{\mathcal{D}_k^q(F_\bu)} \, . 
$$   

In what follows, we will call $S^{k,q}_{F_{u}}$ and $S^{k,q}_{F_{T_u}}$ the first-order and total kernel-based sensitivity indices, respectively. Indeed, we are going to see that $S^{k,q}_{F_{u}}$ and $S^{k,q}_{F_{T_u}}$ represent some well-known first-order and total sensitivity indices for some kernels (see Section \ref{sec:link}). Formal properties of such dependent measures are given in Corollary~\ref{coro:proimnor}.

\begin{corollary} \label{coro:proimnor} 
Let $k$ be an IMK given in Lemma \ref{lem:csak}. Assume that (A1)-(A2) hold. Then,\\   
       
$\quad$ (i) $S^{k,q}_F  \in [0, \, 1], \; \forall \, F \in \mathcal{F}$; \\ 
 
$\quad$ (ii) $S^{k,q}_{F_{T_u}} =0$  iff $\M(\bo{\X})$ is independent of $\bo{\X}_u$; \\ 
  
$\quad$ (iii) $S^{k,q}_{F_{T_u}} =1$ iff $\M(\bo{\X}) = \M(\bo{\X}_u)$; \\

$\quad$ (iv) $S^{k,q}_{F_u}  \leq S^{k,q}_{F_{T_u}}$. 
\end{corollary}   
\begin{preuve} 
See Appendix \ref{app:coro:proimnor}. 
\begin{flushright}  
$\Box$   
\end{flushright} 
\end{preuve}
         
In view of Corollary \ref{coro:proimnor}, an equivalent null hypothesis for the independence test  between the inputs $\bo{\X}_u$ and the outputs is given by
$$
 H_0^{''}: \quad  S^{k,q}_{F_{T_u}} = 0 \, ,   
$$
and the associated test statistic will rely on the estimator of $ S^{k,q}_{F_{T_u}}$. Indeed, performing a statistical test of independence requires an empirical statistic and its distribution under the null hypothesis. 
              
\section{Empirical test statistic and dependent measures}\label{sec:empkbsi}
This section aims at providing empirical dependent measures, including empirical kernel-based sensitivity indices (Kb-SIs) and the test statistic. Note that the first-order AF $\bo{G}_u^{fo}$ and the total AF $\bo{G}_u^{tot}$ for all $u \subseteq \{1, \ldots, d\}$ will lead to the first-order and total Kb-SIs. For concise notations and when there is no ambiguity, we are going to use $S_{u}^{k,q}$, $S_{T_u}^{k,q}$ instead of $S_{F_u}^{k,q}$, $S_{F_{T_u}}^{k,q}$.\\       
         
For computing the dependent measures defined in Section  \ref{sec;impmea}, we are given two i.i.d. samples, that is,  $ \left\{ \left(\bo{\X}_{i,u}, \bo{Z}_{i}, \bo{\X}_{i,u}', \bo{Z}_{i}'  \right) \right\}_{i=1}^{m_1}$ and 
 $ \left\{ \left(\bo{\X}_{i,u}, \bo{Z}_{i}, \bo{\X}_{i,u}', \bo{Z}_{i}'  \right) \right\}_{i=1}^m$  from the random vector $\left(\bo{\X}_{u}, \bo{Z}, \bo{\X}_{u}', \bo{Z}' \right)$, where the four components are mutually independent. Define      
$$
\mu_{k}^{fo} := \esp\left[k\left(\MG^{fo}_u(\bo{\X}_{u}),\, \MG^{fo}_u(\bo{\X}_{u}') \right) \right] ; 
\qquad 
\mu_{k}^{tot} := \esp\left[k\left(\MG^{tot}_u(\bo{\X}_{u} , \bo{Z}),\, \MG^{tot}_u(\bo{\X}_{u}' , \bo{Z}') \right) \right]  \, ;   
$$  
$$          
\sigma_{k}^{fo} := \var\left[k\left(\MG^{fo}_u(\bo{\X}_{u}),\, \MG^{fo}_u(\bo{\X}_{u}') \right) \right] ; \,   
\qquad \,  
\sigma_{k}^{tot} := \var\left[k\left(\MG^{tot}_u(\bo{\X}_{u} , \bo{Z}),\, \MG^{tot}_u(\bo{\X}_{u}' , \bo{Z}') \right) \right] \, ;
$$      
$$
\mu_k^c :=  \esp\left[k\left(\MG(\bo{\X}_{u} , \bo{Z}) - \esp\left[\MG(\bo{\X}_u , \bo{Z}) \right],\, \MG(\bo{\X}_{u}' , \bo{Z}') - \esp\left[\MG(\bo{\X}_u , \bo{Z}) \right] \right) \right] \, .    
$$      
Also, recall that the law of large numbers (LLN) ensures the convergence in probability of the following estimators when $m_1 \to \infty$:        
$$
\widehat{\mu}(\bo{Z}) = \frac{1}{m_1} \sum_{i=1}^{m_1} \MG(\bo{\X}_{i,u} , \bo{Z}) \, \xrightarrow{P} \,  
\esp_{\bo{\X}_u}\left[\MG(\bo{\X}_u , \bo{Z}) \right];  
\quad      
\widehat{\mu}  := \frac{1}{m_1} \sum_{i=1}^{m_1} \MG(\bo{\X}_{i, u} , \bo{Z}_i) \,  \xrightarrow{P} \,
\esp\left[\MG(\bo{\X}_u , \bo{Z}) \right] \, ; 
$$   
$$  
\widehat{\mu}(\bo{\X}_u) := \frac{1}{m_1} \sum_{i=1}^{m_1} \MG(\bo{\X}_{u} , \bo{Z}_i) \, \xrightarrow{P} \,  
\esp_{\bo{Z}}\left[\MG(\bo{\X}_u , \bo{Z}) \right]  \, .  
$$    
  
 Using the plug-in approach, we provide the estimators of $\mu_k^{fo}, \; \mu_k^{tot}$, $\sigma^{tot}_k$ and their statistical properties in Theorem \ref{theo:teststat}. 
\begin{theorem} \label{theo:teststat} 
Let  $k \in \mathcal{K}_{IM}$ be a differentiable kernel almost everywhere (A3), and assume (A1)-(A2) hold. If $ m_1, m \to \infty$, then \\  
    
$\quad$ (i) a consistent estimator of $\mu_k^{tot}$ is given by 
\begin{equation}     
\widehat{\mu_k^{tot}} :=  \frac{1}{m} \sum_{i=1}^m k\left(\MG(\bo{\X}_{i,u} , \bo{Z}_i) - \widehat{\mu}(\bo{Z}_i),\, \MG(\bo{\X}_{i,u}' , \bo{Z}_i') - \widehat{\mu}(\bo{Z}_i') \right)  \, ; 
\end{equation}    
$$
\sqrt{m}\left(\widehat{\mu_k^{tot}} - \mu_k^{tot}\right) \, \xrightarrow{D} \, \mathcal{N} \left(0, \sigma_k^{tot} \right) \, .
$$    

$\quad$ (ii) A consistent estimator of $\sigma_k^{tot}$ is given by 
\begin{equation} \label{eq:estsigtotk}  
\widehat{\sigma_{k}^{tot}} := \frac{1}{m-1} \sum_{i=1}^m\left[ k\left(\MG(\bo{\X}_{i,u} , \bo{Z}_i) - \widehat{\mu}(\bo{Z}_i),\, \MG(\bo{\X}_{i,u}' , \bo{Z}_i') - \widehat{\mu}(\bo{Z}_i') \right) 
- \widehat{\mu_k^{tot}} \right]^2  \, ,      
\end{equation}  
   
$\quad$ (iii) A consistent estimator of $\sigma_k^{tot}$ under the null hypothesis is given by 
\begin{equation} \label{eq:estsigtotk}  
\widehat{\sigma_{k, H_0}^{tot}} := \frac{1}{m} \sum_{i=1}^m\left[ k\left(\MG(\bo{\X}_{i,u} , \bo{Z}_i) - \widehat{\mu}(\bo{Z}_i),\, \MG(\bo{\X}_{i,u}' , \bo{Z}_i') - \widehat{\mu}(\bo{Z}_i') \right) 
- k(\bo{0}, \bo{0}) \right]^2  \, .  
\end{equation} 
 
$\quad$ (iv) A consistent estimator of $\mu_k^{fo}$ is given by 
\begin{equation}     
\widehat{\mu_k^{fo}} :=  \frac{1}{m} \sum_{i=1}^m k\left(\widehat{\mu}(\bo{\X}_{i,u}) - \widehat{\mu},\, \widehat{\mu}(\bo{\X}_{i,u}') - \widehat{\mu} \right)  \, ; 
\qquad
\sqrt{m}\left(\widehat{\mu_k^{fo}} - \mu_k^{fo}\right) \, \xrightarrow{D} \, \mathcal{N} \left(0, \sigma_k^{fo} \right) \, .   
\end{equation}     
\end{theorem}            
\begin{preuve} 
See Appendix \ref{app:theo:teststat}. 
\begin{flushright}
$\Box$  
\end{flushright} 
\end{preuve}   
      
Based on Theorem \ref{theo:teststat}, we derive i) the estimators of Kb-SIs in Corollary \ref{coro:estksi}, and ii) the empirical test statistic under the null hypothesis and its asymptotic distribution in Corollary \ref{coro:teststat2}.
 Estimating the Kb-SIs for at least the $d$ input variables $\X_j$ with $j=1, \ldots, d$ or every subset of inputs will require different samples of the form $\left(\bo{\X}_{u}, \bo{Z}, \bo{\X}_{u}', \bo{Z}' \right)$. We then use $M\geq m$ for the size of the sample that can be used to estimate $\mu_k^c$, that is,  
$$  
\widehat{\mu_k^c} := \frac{1}{M} \sum_{i=1}^M k\left(\MG(\bo{\X}_{i,u} , \bo{Z}_i) - \widehat{\mu},\, \MG(\bo{\X}_{i,u}' , \bo{Z}_i') - \widehat{\mu} \right) \, \xrightarrow{P} \, \mu_k^c \, . 
$$ 
         
\begin{corollary} \label{coro:estksi}  
Let $N \sim \mathcal{N}\left(0, 1 \right)$ be a Gaussian variable, and assume (A1)-(A3) hold.  \\  
  
$\quad$ (i) The consistent estimators of $S_{u}^{k,q}$ and $S_{T_u}^{k,q}$ are given by 
\begin{equation}  \label{eq:estksifo} 
\widehat{S_{u}^{k,q}} :=   \left|\frac{\frac{1}{m} \sum_{i=1}^m k\left(\widehat{\mu}(\bo{\X}_{i,u}) - \widehat{\mu},\, \widehat{\mu}(\bo{\X}_{i,u}') - \widehat{\mu} \right) - k(\bo{0}, \bo{0})}{ \frac{1}{M} \sum_{i=1}^M k\left(\MG(\bo{\X}_{i,u} , \bo{Z}_i) - \widehat{\mu},\, \MG(\bo{\X}_{i,u}' , \bo{Z}_i') - \widehat{\mu} \right) -k(\bo{0}, \bo{0})} \right|^q \, , 
\end{equation}  
\begin{equation} \label{eq:estksitot}  
\widehat{S_{T_u}^{k,q}} := \left|\frac{\widehat{\mu_k^{tot}} -k(\bo{0}, \bo{0})}{\frac{1}{M} \sum_{i=1}^M k\left(\MG(\bo{\X}_{i,u} , \bo{Z}_i) - \widehat{\mu},\, \MG(\bo{\X}_{i,u}' , \bo{Z}_i') - \widehat{\mu} \right) -k(\bo{0}, \bo{0})} \right|^q \, . 
\end{equation}  
   
$\quad$ (ii) If $q=1$, $m_1, m, M \to \infty$ with $\frac{m}{M} \to 0;\, \frac{m_1}{M} \to 0$, then we have the following asymptotic distributions:
$$  
\sqrt{m}\left(\widehat{S_{u}^{k,1}} - S_{u}^{k,1} \right) \, \xrightarrow{D} \, \mathcal{N} \left(0, \frac{\sigma^{fo}_k}{\left(\mu_k^c- k(\bo{0}, \bo{0})\right)^2} \right)\, ; 
$$   
$$ 
\sqrt{m}\left(\widehat{S_{T_u}^{k,1}} - S_{T_u}^{k,1} \right) \, \xrightarrow{D} \,  \mathcal{N} \left(0, \frac{\sigma^{tot}_k}{\left(\mu_k^c- k(\bo{0}, \bo{0})\right)^2} \right) \, .      
$$      
\end{corollary} 
\begin{preuve} 
See Appendix \ref{app:coro:estksi}. 
\begin{flushright}  
$\Box$     
\end{flushright} 
\end{preuve}              
       
\begin{corollary}  \label{coro:teststat2} 
Let $N \sim \mathcal{N}\left(0, 1 \right)$ be a Gaussian variable, and assume (A1)-(A3) hold. \\  
   
$\quad$ (i) If $m_1, m \to  \infty$, then a test statistic under the null hypothesis is given by 
\begin{equation}   \label{eq:teststat2} 
T_{m, H_0}^{q '} := \left|\frac{\widehat{\mu_k^{tot}} - k(\bo{0}, \bo{0})}{\sqrt{\widehat{\sigma_{k, H_0}^{tot}}/m}}\right|^q \, \xrightarrow{D}\, |N|^q \, ,    
\end{equation} 

$\quad$ (ii) If $m_1, m , M \to  \infty$ with $\frac{m_1}{M}, \frac{m_1}{M} \to 0$, then a discrepancy-based test statistic under the null hypothesis is given by 
\begin{equation}   \label{eq:teststat3} 
T_{m, H_0}^{q}  :=   m^{^{\frac{q}{2}}}   \left( \frac{\widehat{\sigma_{k, H_0}^{tot}}}{\left(\widehat{\mu_k^c}- k(\bo{0}, \bo{0})\right)^2} \right)^{-\frac{q}{2}} \, \widehat{S_{T_u}^{k, q}} \; 
\xrightarrow{D}  \;   |N|^q \, . 
\end{equation}   
\end{corollary}        
Using Theorem \ref{theo:teststat} and Corollary \ref{coro:estksi}, the results provided in Corollary~\ref{coro:teststat2} are straightforward by applying the Slutsky theorem. For instance, we can see that 
$$
m^{^{\frac{q}{2}}} \, \widehat{S_{T_u}^{k, q}} \; 
\xrightarrow{D}  \left( \frac{\sigma^{tot}_k}{\left(\mu_k^c- k(\bo{0}, \bo{0})\right)^2} \right)^{q/2} \;   |N|^q \, , 
$$ 
under the null hypothesis. We will rely on $T_{m, H_0}^{q}$ given by (\ref{eq:teststat3}) for performing independence tests, as it is linked to the total Kb-SIs. Thus, the estimators $\widehat{\sigma_{k, H_0}^{tot}}$, $\widehat{\mu_k^c}$ and $\widehat{S_{T_u}^{k, q}}$ are going to be used for performing such independence tests. The testing procedure consists in computing $T_{m, H_0}^{q}$, and then comparing the value obtained to the critical value of $ |N|^q$ at a given threshold such as $\alpha=5\%$. In general, this critical value is the empirical quantile  of the distribution of $|N|^q$ associeted with $1-\alpha$. For $q=1$ or $q=2$, the quantile of the chi or chi-squared distribution of degree $1$ can be directly used.    
               
\section{Links with other importance measures}  \label{sec:link}

For a $\NN$-dimensional random vectors $\bo{G} \sim F$ and the kernels of the form $k (\bo{r}, \bo{r}') := \phi(\bo{r})^\T \phi(\bo{r}')$, the discrepancy measure  
$\mathcal{D}_{k}^{1/2}(F)$ becomes
$$ 
\mathcal{D}_{k}^{1/2}(F) = \esp_{F}\left[\phi(\bo{G})\right] \, .  
$$ 
          
\subsection{Variance-based importance measure} 
For  the kernel $k_{l_2^2} (\bo{r}, \bo{r}') := \norme{\bo{r}}^2 \norme{\bo{r}'}^2$, the associated discrepancy measure   
$ 
\mathcal{D}_{k_{l_2^2}}^{1/2}(F) = \esp_{F}\left[\norme{\bo{G}}^2 \right]
$ 
leads to the generalized sensitivity indices (GSIs) of the first-type, including Sobol' indices (see \cite{sobol93,saltelli00,lamboni11,gamboa14} for independent variables and \cite{lamboni21,lamboni21ar} for dependent and correlated variables). Likewise, the kernel $k_2 (\bo{r}, \bo{r}') := \left<\bo{r},\,  \bo{r}'\right>^2_{\R^\NN}$ and the associated measure   
$
\mathcal{D}_{k_2}^{1/2}(F) = \esp_{F}^{1/2}\left[k_2 (\bo{G}, \bo{G}')\right]
$    
lead to the second-type GSIs (see \cite{lamboni19,lamboni21,lamboni21ar}).

\begin{rem}          
For any AF $\bo{G} \sim F$, it is worth noting that \\     
   
 (i) $\mathcal{D}_{k_{l_2^2}}^{1/2}(F)$ is the first-order Taylor approximation of $\mathcal{D}_{k_g}(F)$ with  $k_g(\bo{r}, \bo{r}') :=e^{-0.5 \norme{\bo{r}- \bo{r}'}^2}$ the Gaussian kernel;\\   
    
 (ii) $\mathcal{D}_{k_2}^{1/2}(F)$ is the second-order Taylor approximation of $\mathcal{D}_{k_e}(F)$ with $k_e(\bo{r}, \bo{r}') :=e^{\sqrt{2}\left<\bo{r},\,  \bo{r}'\right>}$ the exponential kernel.  
\end{rem}    
 
Thus, the well-known variance-based SIs are the approximations of the kernel-based SIs associated with some characteristic IMks.  
    
\subsection{Maximum mean discrepancy and Energy distance} \label{sec;linkmmd}
For any AF $\bo{G} \sim F$ and the centered kernel at zero given by
$
\bar{k}(\bo{r}, \bo{r}') = k(\bo{r}, \bo{r}') - k(\bo{0}, \bo{r}') - k(\bo{r}, \bo{0}) + k(\bo{0}, \bo{0}) 
$; 
we can see that       
\begin{equation}
\mathcal{D}_{\bar{k}}(F) =  \esp_F\left[k(\bo{G}, \bo{G}')\right] -2 \esp_F\left[k(\bo{G}, \bo{0})\right] 
+ k(\bo{0}, \bo{0}) = \normh{\esp_F\left[k(\cdot, \bo{G})\right] - k(\cdot, \bo{0})}^2  \, , \nonumber 
\end{equation}
 is the squared MMD between the distribution $F$ and the Dirac measure $\delta_{\bo{0}}$.  \\ 

Kernels induced by the semimetric $\norme{\bo{r}- \bo{r}'}^\alpha$ with $0 < \alpha < 2$, that is, 
$$  
k_d^\alpha (\bo{r}, \bo{r}') := \frac{1}{2} \left(\norme{\bo{r}}^\alpha + \norme{\bo{r}'}^\alpha - \norme{\bo{r}-\bo{r}'}^\alpha \right) 
$$     
are centered at zero and belong to the set of kernels that guarantee the independence criterion. For such kernels, the measure
$
\mathcal{D}_{k_d^1}(F) = \esp_{F}\left[k_d^1 (\bo{G}, \bo{G}') \right]
$      
is twice the squared energy distance between the distribution $F$ and the Dirac delta measure $\delta_{\bo{0}}$ (see \cite{rizzo16} for more details). Recall that such kernels must be convex or concave on the support of AFs to be IMKs.\\      
      
While the MMD and distance energy are used for testing independence between random vectors, it comes out that additional conditions on the associated kernels are needed in order to obtain ANOVA-compatible kernels and importance dependent measures between the inputs and the outputs.   
    
\section{Importance measure kernels  based on log-concave kernels} \label{app:conkerp}
This section shows how one can control log-concave kernels of the form $\exp\left(-\alpha \psi(\bo{y}, \bo{y}') \right)$ to obtain concave kernels, where $\alpha >0$ and $\psi$ is a convex function. Note that some radial-based characteristic kernels such as the Laplacian and Gaussien kernels are log-concave kernels. To that end, we use $\partial \psi$ for the subgradient of $\psi$ (\cite{boyd04}) and $H_\psi$ for the Hessian of $\psi$. When $\psi$ is differentiable, the subgradient $\partial \psi$ is equal to the gradient $\nabla \psi$.    
  
\begin{lemma} \label{lem:conker}     
Let   $\mathcal{\X}$ be the support of the outputs $\bo{Y}$  and $0 <\epsilon \leq 1$.  Assume $\psi$ is convex and continuous. \\       
    
$\quad$ (i) The function $\bo{y} \mapsto \exp\left(-\alpha \psi(\bo{y}, \bo{y}') \right)$ is concave on $\mathcal{\X}$ whenever         
$$
\alpha \leq  
\inf_{\bo{y} \in \mathcal{\X}}   
\inf_{\bo{y}' \in \mathcal{\X}}
 \inf_{\bo{z} \in \mathcal{\X}} \frac{\epsilon}{\norme{\bo{y}-\bo{y}'} \norme{\partial \psi(\bo{y}', \bo{z})}} \, .
$$
 
$\quad$ (ii) If $\psi$ is twice differentiable, then $\exp\left(-\alpha \psi(\bo{y}, \bo{y}') \right)$ is concave on $\mathcal{\X}$ when 
$$
\alpha \leq \max\left( 
\inf_{\bo{y} \in \mathcal{\X}}   
\inf_{\bo{y}' \in \mathcal{\X}}
 \inf_{\bo{z} \in \mathcal{\X}} \frac{\epsilon}{\norme{\bo{y}-\bo{y}'} \norme{\partial \psi(\bo{y}', \bo{z})}}
 \, ,  \,   
\inf_{\bo{y} \in \mathcal{\X}} 
\inf_{\bo{y}' \in \mathcal{\X}}
  \inf_{\bo{b} \in \mathcal{\X}} \frac{\bo{b}^\T  H_\psi(\bo{y}, \bo{y}') \bo{b} }{\left( \bo{b}^\T \nabla \psi(\bo{y}, \bo{y}')\right)^2}
	\right) \, .  
$$
\end{lemma}    
\begin{preuve}
See Appendix \ref{app:lem:conker}. 
\begin{flushright}
$\Box$
\end{flushright}
\end{preuve} 
   
For a bounded support $\mathcal{\X}$, that is, $\norme{\bo{y}} \leq C$ for all $\bo{y} \in \mathcal{\X}$ with $C>0$ a constant, the first condition becomes 
$$  
\inf_{\bo{y} \in \mathcal{\X}} 
\inf_{\bo{y}' \in \mathcal{\X}}
 \frac{\epsilon}{2C \norme{\partial \psi(\bo{y}, \bo{y}')}} >0 \, . 
$$      
For instance,  the Laplacian kernel given by $e^{-\alpha \norml{\bo{y}- \bo{y}'}}$ is concave on the bounded support $\mathcal{\X}$ for every $\alpha$ satisfying 
$
\alpha \leq \frac{\epsilon}{2C \sqrt{\NN}}   
$. 
Likewise, the Gaussian kernel given by $e^{-\alpha \norme{\bo{y}- \bo{y}'}^2}$ is concave on $\mathcal{\X}$ when
$
\alpha \leq \frac{\epsilon}{8 C^{2}}
$. \\ 
 
For an unbounded support $\mathcal{\X}$, which is the case of most log-concave probability measures, the parameter $\alpha$ can be chosen so that the inequalities provided in Lemma \ref{lem:conker} fail with a small value of probability, that is, about $0.05$ (see Corollary \ref{coro:conkerunb}).  
\begin{corollary} \label{coro:conkerunb}   
Let $\bo{Y}', \bo{Y}''$ be two i.i.d. copies of the outputs $\bo{Y}$;  $0 <\tau \leq 5\%$ and $0< \epsilon \leq 1$. Assume that $\psi$ is convex and continuous. \\  

 $\quad$ (i) The kernel $k$ is concave on $\mathcal{\X}$ with high probability ($\geq 1-\tau$) when
$$  
\alpha  \leq \max\left(\frac{\tau}{\epsilon\,\esp\left[\norme{\bo{Y}-\bo{Y}'} \norme{\partial \psi(\bo{Y}', \bo{Y}'')} \right]},\, 
\frac{\sqrt{\tau}}{\epsilon \sqrt{ \esp\left[\norme{\bo{Y}-\bo{Y}'}^2 \norme{\partial \psi(\bo{Y}', \bo{Y}'')}^2 \right]}} \right) \, .   
$$    
   
$\quad$ (ii)  If $\psi$ is twice differentiable, then $k$ is concave with high probability when 
$$ 
\alpha \leq \max \left( \frac{\tau}{\esp\left[\frac{\left(\nabla \psi(\bo{Y}, \bo{Y}')^\T \, \bo{Y}''\right)^2}{\bo{Y}^{'' \,\T}  H_\psi(\bo{Y}, \bo{Y}') \bo{Y}''} \right]}, 
\,
\frac{1}{1 -\tau} \esp\left[ \frac{\bo{Y}^{'' \,\T}  H_\psi(\bo{Y}, \bo{Y}') \bo{Y}''}{\left(\nabla \psi(\bo{Y}, \bo{Y}')^\T \, \bo{Y}''\right)^2} \right] \right)  \, . 
$$
\end{corollary}        
\begin{preuve}
See Appendix \ref{app:coro:conkerunb}. 
\begin{flushright} 
$\Box$  
\end{flushright}
\end{preuve}
 
Thus, the Laplacian kernel is concave on the support $\mathcal{\X}$ with higher probability ($>1-\tau$) when 
$
\alpha \leq \frac{\sqrt{\tau}}{\epsilon \sqrt{2\NN \trace\left(\var[\bo{Y}] \right)}}    
$. 
In the case of the Gaussian kernel, the condition becomes 
$
\alpha \leq \frac{\tau}{4\epsilon \trace\left(\var[\bo{Y}] \right)}       
$. For a given $\tau$, We may choose $\epsilon=\sqrt{\tau}$ for the first condition and $\epsilon=\tau$ for the second one. 
           
\section{Simulation study}   \label{sec:test}
To illustrate our approach, we consider two functions and the following kernels: the quadratic kernel $k_2(\bo{r}, \bo{r}') =\left(\bo{r}^\T \bo{r}'\right)^2$, the $L_1$-based kernel $k_{l_1}(\bo{r}, \bo{r}') =\left|\left|\bo{r} \right|\right|_1 \left|\left|\bo{r}'\right|\right|_1$, the Gaussian kernel $k_G(\bo{r}, \bo{r}') =\exp\left(-\alpha_1 \, \norme{\bo{r} - \bo{r}'}^2\right)$ and the Laplacian kernel $k_L(\bo{r}, \bo{r}') =\exp\left(-\beta \, \left|\left| \bo{r} - \bo{r}'\right|\right|_1 \right)$. In this section, $\alpha_1$ and  $\beta$ were chosen according to Corollary \ref{coro:conkerunb} using the variance of the model output(s).  
            
\subsection{Sobol' function} \label{sec:test1} 
Consider a model that includes ten independent variables following the uniform distribution, that is, $\X_j \sim \mathcal{U}(0,\, 1)$ with $j=1, \ldots, d$, and given by   
$$ 
\M(\bo{\X}) = \prod_{j=1}^{d=10}\frac{ |4\, \x_j \,- \,2| \,+ \, \bo{a}[j]}{1 \,+\, \bo{a}[j]} \, ,   \quad  \mbox{with} \;   \bo{a} :=[0,\, 0,\, 6.52, \ldots, 6.52]^\T  \, .
$$
The variance of $\M(\bo{\X})$ is $\var\left[\M(\bo{\X}) \right] = 0.863$, and the Kb-SIs were computed using $m_1=1000$, $M=m=2000$, $q=1/2$, $\alpha_1 =1/8 < 0.29$ and $\beta=1/4 <  0.76$ (see Corollary \ref{coro:conkerunb}). The estimated Kb-SIs are reported in Table \ref{tab:test1}, including the Kb-SIs associated with the exponential kernel, that is, $k_e(\bo{r}, \bo{r}') :=e^{\sqrt{2}\left<\bo{r},\,  \bo{r}'\right>}$. \\        
              
The independence statistical tests based on (\ref{eq:teststat3})) reveal that the output $\M(\bo{\X})$ depends on all the input variables for most of the kernels, except the exponential kernel. According to the statistical test based on the exponential kernel, $\M(\bo{\X})$ depends only on $\X_1$ and $\X_2$. For other kernels, it comes out that we have the same ranking of inputs using the total Kb-SIs. For selecting the  most influential variables, it is common to fix a threshold $T$. When $T=0.1$, it appears that all the inputs are important according the the Gaussian, Laplacian and somehow the L1-based kernels. Sobol' indices (quadratic kernel) and the exponential Kb-SIs identify $\X_1$ and $\X_2$ as the most important variables. Such differences are due to the fact that different kernels capture different information of AFs. For instance,  it is known that small norms such as the the L1-based kernel capture slow variations of AFs.   
   
\begin{table}[http]     
\centering 
\begin{tabular}{lcccccccccc}
  \hline 
	  \hline 
Kernels & X1 & X2 & X3 & X4 & X5 & X6 & X7 & X8 & X9 & X10  \\ 
  \hline   
\multicolumn{11}{c}{First-order Kb-SIs} \\
 L1-based & 0.678 & 0.679 & 0.090 & 0.090 & 0.090 & 0.090 & 0.090 & 0.090 & 0.090 & 0.090 \\ 
  Quadratic & 0.393 & 0.393 & 0.007 & 0.007 & 0.007 & 0.007 & 0.007 & 0.007 & 0.007 & 0.007 \\ 
  Gaussian & 0.695 & 0.694 & 0.097 & 0.098 & 0.097 & 0.097 & 0.097 & 0.097 & 0.097 & 0.097 \\ 
  Laplacian & 0.854 & 0.853 & 0.328 & 0.329 & 0.327 & 0.327 & 0.327 & 0.328 & 0.327 & 0.327 \\ 
  Exponential & 0.106 & 0.107 & 0.001 & 0.003 & 0.002 & 0.002 & 0.002 & 0.001 & 0.001 & 0.001 \\ 
\hline    
\multicolumn{11}{c}{Total Kb-SIs} \\ 	
  L1-based & 0.675 & 0.678 & 0.089 & 0.089 & 0.089 & 0.089 & 0.089 & 0.089 & 0.088 & 0.091 \\ 
  Quadratic & 0.531 & 0.548 & 0.013 & 0.012 & 0.012 & 0.012 & 0.012 & 0.012 & 0.012 & 0.013 \\ 
  Gaussian & 0.787 & 0.787 & 0.131 & 0.131 & 0.131 & 0.130 & 0.130 & 0.131 & 0.131 & 0.131 \\ 
  Laplacian & 0.894 & 0.890 & 0.356 & 0.357 & 0.356 & 0.356 & 0.355 & 0.357 & 0.357 & 0.356 \\ 
   Exponential & 0.174 & 0.188 & 0.002 & 0.002 & 0.002 & 0.005 & 0.004 & 0.003 & 0.002 & 0.004 \\ 
\hline    
\multicolumn{11}{c}{Test statistics values ($T_{m, H_0}^{q}$, (\ref{eq:teststat3}))} \\ 	
 L1-based & 4.966 & 4.907 & 4.236 & 4.362 & 4.355 & 4.314 & 4.285 & 4.289 & 4.374 & 4.271 \\ 
 Quadratic & 3.659 & 3.494 & 2.640 & 2.770 & 2.928 & 2.677 & 2.654 & 2.745 & 2.955 & 2.836 \\ 
 Gaussian & 5.244 & 5.189 & 4.517 & 4.554 & 4.571 & 4.609 & 4.543 & 4.504 & 4.590 & 4.524 \\ 
 Laplacian & 6.034 & 6.008 & 5.688 & 5.699 & 5.695 & 5.710 & 5.702 & 5.707 & 5.711 & 5.692 \\ 
   Exponential & 2.621 & 2.259 & 0.314 & 0.364 & 0.275 & 0.959 & 0.687 & 0.473 & 0.305 & 0.624 \\ 
\hline    
\multicolumn{11}{c}{Critical values} \\ 
 L1-based & 1.393 & 1.412 & 1.379 & 1.386 & 1.394 & 1.387 & 1.408 & 1.393 & 1.398 & 1.370 \\ 
 Quadratic & 1.412 & 1.426 & 1.410 & 1.380 & 1.409 & 1.363 & 1.403 & 1.422 & 1.386 & 1.396 \\ 
 Gaussian & 1.399 & 1.399 & 1.406 & 1.399 & 1.373 & 1.391 & 1.399 & 1.374 & 1.392 & 1.405 \\ 
 Laplacian & 1.387 & 1.409 & 1.404 & 1.394 & 1.401 & 1.409 & 1.395 & 1.411 & 1.393 & 1.409 \\ 
  Exponential & 1.420 & 1.422 & 1.385 & 1.415 & 1.383 & 1.388 & 1.389 & 1.413 & 1.396 & 1.385 \\
\hline    
\multicolumn{11}{c}{Decision about dependence} \\  
 L1-based & Yes & Yes & Yes & Yes & Yes & Yes & Yes & Yes & Yes & Yes \\ 
 Quadratic & Yes & Yes & Yes & Yes & Yes & Yes & Yes & Yes & Yes & Yes \\ 
 Gaussian & Yes & Yes & Yes & Yes & Yes & Yes & Yes & Yes & Yes & Yes \\ 
 Laplacian & Yes & Yes & Yes & Yes & Yes & Yes & Yes & Yes & Yes & Yes \\ 
Exponential & Yes & Yes & No & No & No & No & No & No & No & No \\ 
  \hline   
	\hline         
\end{tabular} 
\caption{Kernel-based sensitivity Indices, values of test statistics and critical values at the threshold $\alpha = 0.05$.}
\label{tab:test1}       
\end{table}

\subsection{Vector-valued function} \label{sec:test2}
\begin{figure}[!hbp]   
\begin{center}
\includegraphics[height=10cm,width=14cm,angle=0]{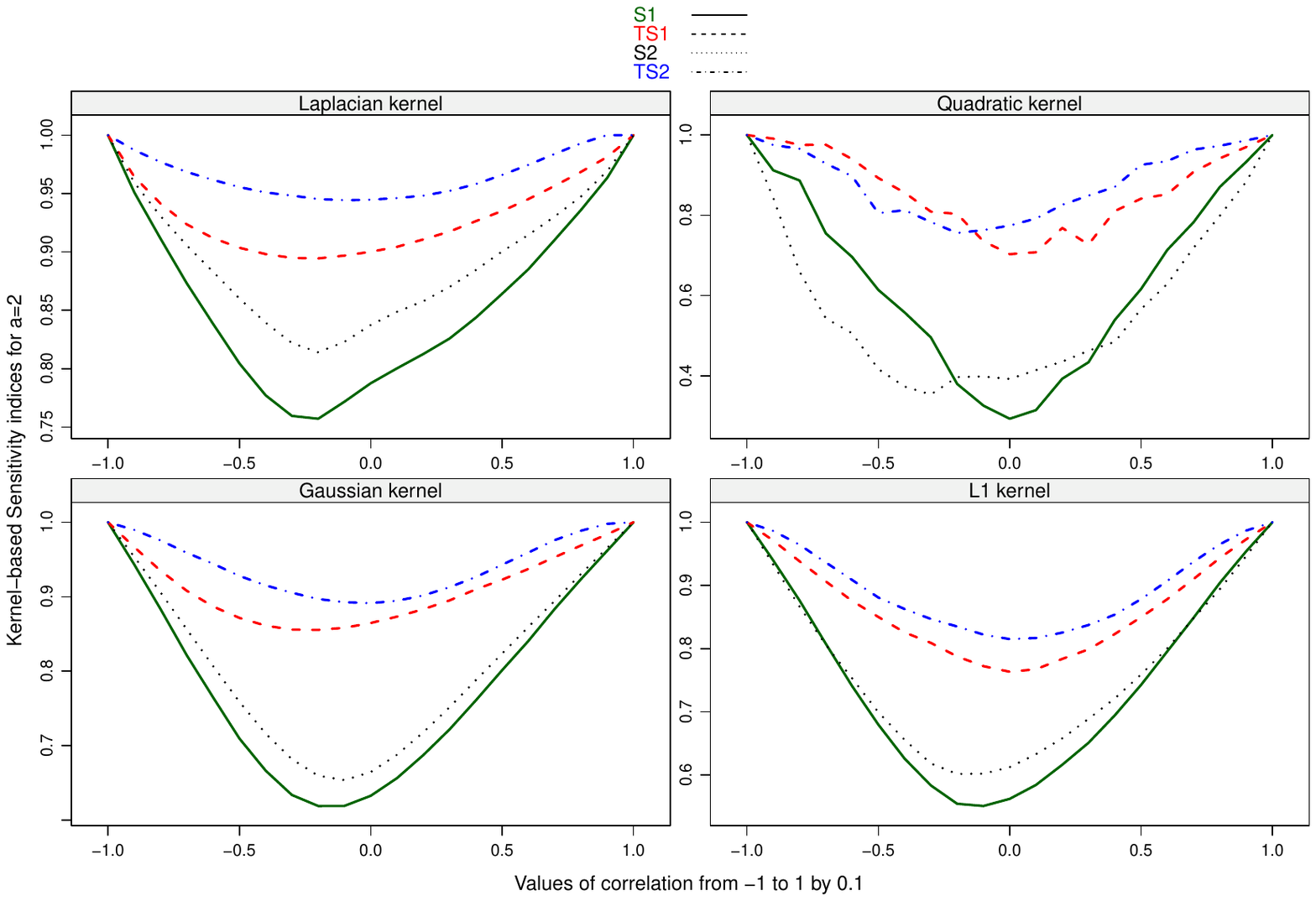}
\end{center} 
\caption{Kernel-based sensitivity indices for q=1/2 (i.e., $S^{k,1/2}$) and a=2}  
 \label{fig:ksi}     
\end{figure}

Consider the following model   
\begin{equation}            
\M(\X_1,\, \X_2) := \left[
\X_1 +\X_2 + a \X_1 \X_2,
\quad  
 \X_1^2 + \sqrt{2}\X_2       
\right]^\T \, ,  
\end{equation}  
which includes two correlated variables, that is, $\X_j \sim \mathcal{N}(0,\, 1),\, j =1,\, 2$ with $\rho$ the correlation coefficient and $a \in \R$. Using the dependency models (\cite{lamboni21,lamboni21ar}), that is, $\X_2 \stackrel{d}{=} \rho \X_1 +\sqrt{1-\rho^2} Z_2$ and $\X_1\stackrel{d}{=}\rho \X_2 +\sqrt{1-\rho^2} Z_1$ with $Z_j \sim \mathcal{N}\left(0,\, 1 \right)$ and $j=1, 2$, the equivalent representations of the model are given by   
$$   
\left[  
 (1+\rho)X_1 + a \rho\X_1^2 + \sqrt{1-\rho^2}Z_2 (1 +  a \X_1),
\quad \, 
\X_1^2 + \sqrt{2}\rho \X_1 + \sqrt{2}\sqrt{1-\rho^2} Z_2 \right]^\T \, , 
$$ 
$$
\left[
(1+\rho)X_2 + a \rho\X_2^2 + \sqrt{1-\rho^2} Z_1(1 + a \X_2),
\quad \, 
	\rho^2 \X_2^2	+ (1-\rho^2)Z_1^2 +2\rho \sqrt{1-\rho^2}Z_1 \X_2 + \sqrt{2} X_2
\right]^\T \, . 
$$  
The first-order and total AFs are given by
$$ 
\BG_1^{fo} =  \left[\begin{array}{c} 
(1+\rho)X_1 + a \rho\X_1^2 -a \rho \\
 \X_1^2 + \sqrt{2}\rho \X_1 -1     \\ 
\end{array}
\right];
\quad   
\BG_1^{tot} =  \left[\begin{array}{c} 
(1+\rho)X_1 + a \rho\X_1^2 + a \sqrt{1-\rho^2}Z_2 \X_1 -a \rho  \\  
 \X_1^2 + \sqrt{2}\rho \X_1 -1  \\  
\end{array} \right] \, ; 
$$ 
$$  
\BG_2^{fo} =  \left[\begin{array}{c} 
(1+\rho)X_2 + a \rho\X_2^2 -a \rho  \\
\rho^2 \X_2^2	+ \sqrt{2} X_2 -\rho^2  \\
\end{array} 
\right]; 
\quad   
\BG_2^{tot} =  \left[\begin{array}{c} 
(1+\rho)X_2 + a \rho\X_2^2 + a\sqrt{1-\rho^2} Z_1 \X_2 -a \rho  \\  
\rho^2 \X_2^2	+ 2\rho \sqrt{1-\rho^2}Z_1 \X_2 + \sqrt{2} X_2 -\rho^2   \\ 
\end{array} \right] \, . 
$$     
When $\rho=0$, both inputs are independent, and we can see that the components of $\BG_2^{fo}$ and $\BG_2^{tot}$ are correlated while those of $\BG_1^{fo}$ and $\BG_1^{tot}$ are clearly not, but they are dependent. We can also check that such AFs are not Gaussian-distributed. Figure \ref{fig:ksi} compares estimates of kernel-based SIs associated with the four kernels using $\alpha_1 =\frac{1}{4(6+a^2)}$ and $\beta=\frac{1}{\sqrt{4(6+a^2)}}$ with $6+a^2$ the trace of the covariance of $\M(\X_1,\, \X_2)$ when $\rho=0$ (see \cite{lamboni21}).    
We can see that the total Kb-SIs are always greater than the first-order indices, as expected. While the Gaussian, Laplacian and $L_1$ Kb-SIs give the same ranking of inputs, the quadratic Kb-SIs inversely identify $\X_1$ as the most important input for negative values of the correlation coefficient. 
                             
\section{Conclusion} \label{sec:con}   
In this paper, we have proposed a new dependent measure of association between the model inputs and outputs and the associated statistical test of independence by making use of the total ANOVA functional (AF) of inputs and kernel methods. The proposed test statistic and dependent measures are i) well-suited for any output domain or object and models with non-independent variables, ii) much flexible for explicitly including specific moments of AFs, iii) able to distinguish two different distributions of AFs. Regarding the test statistics, we have provided generic kernels that guarantee the independence criterion between the outputs and the inputs. Most of characteristic kernels are equivalent kernels for the independence criterion. It comes out that convex and some concave kernels that guarantee the independence criterion lead to interesting dependent measures according to R\'enyi' axioms.\\ 
               
In uncertainty quantification, the variance-based SA (Sobol' indices; generalized sensitivity indices (\cite{lamboni11,gamboa14,gamboa14b,lamboni19}) and dependent generalized sensitivity indices (dGSIs) (\cite{lamboni21,lamboni21ar})) and the Owen $L^p$ measure are  special cases of the proposed dependent measures by properly choosing the kernels. Moreover, it comes out that some kernels lead to the application of the maximum mean discrepancy and distance energy  (\cite{borgwardt06,gretton07,gretton12,rizzo16}). For the choice of kernels, we should prefer the importance-measure kernels that detect independence between inputs and outputs as effectively as possible, that is,  account for the necessary and sufficient high-order moments such as the first and second-order moments for Gaussian-distributed AFs. For other distributions of AFs, the Laplacian, Gaussian and exponential kernels that account for all the moments of any distribution can be used since such kernels are importance-measure kernels for some values of their parameters. In general, Sobol' indices and dGSIs are a first-order (resp. second-order) approximation of the proposed dependent measures associated with the Gaussian (resp. exponential) kernel. \\         
   
The computations of AFs require evaluating conditional expectations of the model outputs given some inputs, and such evaluations may be time demanding to converge. In next future, it is interesting to investigate efficient methods for computing conditional expectations in high dimension. Gaussian process emulator (\cite{oakley04,conti10}) of the models or new emulator based on the paper \cite{lamboni22} or the Nadaraya-Watson kernel estimator for high-dimensional nonparametric regression (\cite{conn19}) are going to be investigated. 
                               
  \begin{appendices}      
\section{Proof of Lemma \ref{lem:cripr}} \label{app:lem:cripr}
We can see that the output $\M(\bo{\X})$ is independent of 
$ \bo{\X}_u$ implies  
$$ 
\M(\bo{\X}) \stackrel{d}{=} \M(\bo{\X}_u, r_u(\bo{\X}_u, \bo{Z})) =: 
 \MG(\bo{\X}_u, \bo{Z}) =  \MG(\bo{Z})\, .    
$$          
 Therefore, we have $\MG^{tot}_u(\bo{\X}_u, \bo{Z}) =\MG(\bo{\X}_{u}, \bo{Z}) = \esp_{\bo{\X}_{u}} \left[\MG(\bo{\X}_{u}, \bo{Z}) \right] =\bo{0} \, a.s.$.\\ 
 Conversely, if $\MG^{tot}_u(\bo{\X}_{u}, \bo{Z}) =\bo{0}$, the properties of conditional expectation show that there exists a function $h$ such that 
$$ 
\MG(\bo{\X}_{u}, \bo{Z}) = \esp_{\bo{\X}_{u}} \left[\MG(\bo{\X}_{u}, \bo{Z}) \right] = h(\bo{Z})\, , 
$$       
which means that $\MG(\bo{\X}_{u}, \bo{Z})$ is a function of $\bo{Z}$ only, and the result holds.  

\section{Proof of Lemma \ref{lem:extkerind}} \label{app:lem:extkerind}
Bearing in mind Definition \ref{def:equitest}, we want to show that  
\begin{equation}  
 \esp\left[k\left(\MG^{tot}_u (\bo{\X}_{u}, \bo{Z}), \,  \MG^{tot}_u (\bo{\X}_{u}', \bo{Z}') \right)\right] -k(\bo{0}, \bo{0})  = 0  \, , 
\Longrightarrow \, \MG^{tot}_u (\bo{\X}_{u}, \bo{Z})=\bo{0} \,.
  \nonumber     
\end{equation}       
First, let us start with the kernel $k$. Using the theorem of transfer, we can write 
$$
\esp\left[k\left(\MG^{tot}_u (\bo{\X}_{u}, \bo{Z}), \,  \MG^{tot}_u (\bo{\X}_{u}', \bo{Z}') \right)\right] -k(\bo{0}, \bo{0}) 
=  
\int_{\mathcal{X}^2} k\left(w,  w' \right) \, d(F_{T_u}\otimes F_{T_u} -H \otimes H)(w,w') \, , 
$$        
with $F_{T_u}$ the CDF of $\MG^{tot}_u (\bo{\X}_{u}, \bo{Z})$ and $H$ the CDF of $\delta_{\bo{0}}$. 
For a SPD kernel $k$, when the above identity is zero, it implies that $F_{T_u} \otimes F_{T_u} = H \otimes H$; $F_{T_u} =H$ and $\MG^{tot}_u (\bo{\X}_{u}, \bo{Z})=\bo{0} \, a.s.$.\\  
Second, we are going to use the criterion $\esp_{\bo{G} \sim F, \bo{G}'\sim F}\left[\bar{k}(\bo{G},  \bo{G}')\right] =0 \Rightarrow F= H$. Since $\bar{k}$ is centered at $\bo{0}$, we  have
$$   
\esp\left[\bar{k}\left(\MG^{tot}_u (\bo{\X}_{u}, \bo{Z}), \,  \MG^{tot}_u (\bo{\X}_{u}', \bo{Z}') \right)\right] -\bar{k}(\bo{0}, \bo{0}) 
= 
\esp_{\BG^{tot}\sim F_{T_u}, \BG^{tot '} \sim F_{T_u}}\left[\bar{k}\left(\BG^{tot}, \,  \BG^{tot '} \right)\right] =0 \, ,
$$       
which implies that $F_{T_u} =H$.      

\section{Proof of Lemma \ref{lem:kerequiv}} \label{app:lem:kerequiv}
For Point (i), we are going to use the first criterion of $\mathcal{K}_E$ (see Equation (\ref{eq:setker})). We can write for all $F \in \mathcal{F}$ 
$$    
\esp_{\bo{G} \sim F, \bo{G}' \sim F}\left[\bar{k}\left(\bo{G},  \bo{G}' \right)\right] 
= \normh{\esp_{\bo{G} \sim F}\left[\bar{k}(\cdot, \bo{G}) \right]}^2
=0 \Longleftrightarrow \esp_{\bo{G} \sim F}\left[\bar{k}(\cdot, \bo{G}) \right] = 0\, .
$$      
Note that $\esp_{\bo{G} \sim F}\left[\bar{k}(\cdot, \bo{G}) \right] = 0 = \esp_{W \sim H}\left[\bar{k}(\cdot, W) \right]$ with $H$ the CDF of $\delta_{\{\bo{0}\}}$, and $\bar{k}$ is a characteristic kernel when $k$ is a characteristic one. Thus, Point (i) holds because  $\esp_{\bo{G} \sim F}\left[\bar{k}(\cdot, \bo{G}) \right] = \esp_{W\sim H}\left[\bar{k}(\cdot, W) \right]$ implies $F =H$. \\    
For Point (ii), we are going to use the second criterion of $\mathcal{K}_E$. According to the Bochner Lemma and the Fubini theorem, we can write for independent vectors $Y, Z$
\begin{eqnarray}   
\esp_{Y\sim \mu, Z \sim \mu}\left[k(Y, Z)\right]  &=& \int_{\mathcal{X}^2} k(\bo{y}, \bo{z}) \, d\mu\otimes \mu(\bo{y}, \bo{z}) \nonumber \\ 
& =& \int_{\mathcal{X}^2} \int_{\R^\NN} e^{-i \bo{y}^\T \bo{w}} e^{i \bo{z}^\T \bo{w}} \, d\mu(\bo{y}) d\mu(\bo{z}) d\Lambda(\bo{w}) \nonumber \\ 
&=& \int_{\R^\NN} \left| \int_{\mathcal{X}} e^{-i \bo{y}^\T \bo{w}} \, d\mu(\bo{y}) \right|^2  d\Lambda(\bo{w})\nonumber \, .   
\end{eqnarray}    
Thus, $\esp_{Y\sim \mu, Z \sim \mu}\left[k(Y, Z)\right] =0$ implies that  $\int e^{-i \bo{y}^\T \bo{w}} \, d\mu(\bo{y})=0$ for all $\bo{w} \in Supp(\Lambda)= \R^\NN$. 
For a class of finite and signed Borel measures of the form $\mu(A) := \int_A h(\bo{y}) \, d\bo{y} $ with $h: \R^\NN \to \R$ a measurable function such as a difference of two probability densities, the function  
$$  
 \widehat{h}(\bo{w}) := 
\int e^{-i \bo{y}^\T \bo{w}} \, d\mu(\bo{y}) = \int e^{-i \bo{y}^\T \bo{w}} h(\bo{y}) \,d\bo{y}, \quad \forall\; \bo{w} \in Supp(\Lambda) =\R^\NN
\, , 
$$  
is the Fourier transform of $h(\bo{y})$. As $\widehat{h}(\bo{w})=0$ for all $\bo{w} \in \R^\NN$, we then have $h(\bo{y}) =0$ bearing in mind the inverse Fourier transform. Point (ii) holds because $\mu=0$. 

\section{Proof of Lemma \ref{lem:csak}} \label{app:lem:csak}
Let $\mathcal{H}$ denote an Hilbert space induced by $k$. Without loss of generality, we are going to show the results for $q=1$.\\ 
First, using the convexity of $J(f) :=\normh{f}^2$ with $f \in \mathcal{H}$, we know that there exist a gradient of $\normh{f}^2$ (i.e., $\nabla J(f) := 2f$)  such that for all $f_0 \in \mathcal{H}$  (\cite{boyd04}) 
$$   
\normh{f}^2 - \normh{f_0}^2 \geq \left< 2 f_0,\,  f-f_0 \right>_{\mathcal{H}}\, .
$$  
Second, for $\BG^{tot}_u \stackrel{}{=} \MG^{tot}_u (\bo{\X}_{u}, \bo{Z})$ and $\BG^{fo}_u \stackrel{}{=} \MG^{fo}_u (\bo{\X}_{u})$, we have  (see Equation (\ref{eq:linktotfo}))
$$
\esp_{\bo{Z}}\left[\MG^{tot}_u (\bo{\X}_{u}, \bo{Z}) \right] = \MG^{fo}_u (\bo{\X}_{u})   \, .  
$$  
For Point (i), knowing that for the centered kernel $\bar{k}$, 
$$ 
\mathcal{D}_k(F_{T_u}) = 
\esp\left[\bar{k}\left(\BG^{tot}_u, \,  \BG^{tot '}_u \right)\right] = \normh{\esp\left[k\left(\cdot, \BG^{tot}_u \right)\right]- k\left(\cdot, \bo{0}\right)}^2 \, ,  
$$         
we can write  (bearing in mind that  $k(\bo{0}, \bo{y}') =k(\bo{y}, \bo{0}) =c$) 
\begin{eqnarray}  \label{eq:appdifker}   
& & \mathcal{D}_k(F_{T_u})- \mathcal{D}_k(F_{u}) =  \esp\left[\bar{k}\left(\BG^{tot}_u, \,  \BG^{tot '}_u \right)\right]    
- 
\esp\left[\bar{k}\left(\BG^{fo}_u , \,  \BG^{fo '}_u \right)\right]   \nonumber \\
& \geq &    2\left<\esp\left[k\left(\cdot, \BG^{fo}_u \right)\right]- k\left(\cdot, \bo{0}\right), \;  
 \esp\left[k\left(\cdot, \BG^{tot}_u \right)\right]- \esp\left[k\left(\cdot, \BG^{fo}_u\right)\right] \right>_{\mathcal{H}}    \nonumber \\
	& =&   2  \esp\left[k\left(\BG^{tot}_u, \,  \BG^{fo '}_u\right)\right]
- 2 \esp\left[k\left(\BG^{fo}_u, \,  \BG^{fo '}_u \right)\right] -2 \esp\left[k\left(\BG^{tot}_u, \,  \bo{0}\right)\right] + 2 \esp\left[k\left(\BG^{fo}_u, \,  \bo{0} \right) \right] \nonumber \\ 
	& =&   2  \esp\left[k\left(\BG^{tot}_u, \,  \BG^{fo '}_u\right)\right]
- 2 \esp\left[k\left(\BG^{fo}_u, \,  \BG^{fo '}_u \right)\right]  \nonumber  \, .
\end{eqnarray}         
Point (i) holds using the Jensen inequality and Equation (\ref{eq:linktotfo}). \\  
For (ii), since $\esp [\BG^{tot}_u] =\esp [\BG^{fo}_u]=\bo{0}$ and $k$ is convex, we can write $\mathcal{D}_k(F_{T_u})$ without the absolute symbol thanks to Jensen's theorem, that is,  
$$        
\mathcal{D}_k(F_{T_u}) = 
\esp\left[k\left(\BG^{tot}_u, \,  \BG^{tot '}_u \right)\right] -k(\bo{0}, \bo{0}) = \normh{\esp\left[k\left(\cdot, \BG^{tot}_u \right)\right]}^2 -k(\bo{0}, \bo{0}) \, .
$$         
Using the convexity of $\normh{\cdot}^2$, we can write 
\begin{eqnarray}  \label{eq:appdifker}   
& & \mathcal{D}_k(F_{T_u})- \mathcal{D}_k(F_{u}) 
= \normh{\esp\left[k\left(\cdot, \BG^{tot}_u \right)\right]}^2 - \normh{\esp\left[k\left(\cdot, \BG^{fo}_u \right)\right]}^2  \nonumber \\    
& \geq &    2\left<\esp\left[k\left(\cdot, \BG^{fo}_u\right)\right], \; 
 \esp\left[k\left(\cdot, \BG^{tot}_u \right)\right]- \esp\left[k\left(\cdot, \BG^{fo}_u \right)\right] \right>_{\mathcal{H}}   \nonumber \\
	& =&   2  \esp\left[k\left(\BG^{tot}_u, \,  \BG^{fo '}_u \right)\right]
- 2  \esp\left[k\left(\BG^{fo}_u, \,  \BG^{fo '}_u\right)\right]   \nonumber  \, .       
\end{eqnarray}          
Thus, Point (ii) holds using the Jensen inequality and Equation (\ref{eq:linktotfo}). \\    
For Point (iii), as $k(\bo{0}, \bo{0}) >0$ and $k$ is concave, we have         
$$
\mathcal{D}_k(F_{T_u}) = 
- \esp\left[k\left(\BG^{tot}_u, \,  \BG^{tot '}_u \right)\right] + k(\bo{0}, \bo{0}) = -\normh{\esp\left[k\left(\cdot, \BG^{tot}_u \right)\right]}^2 + k(\bo{0}, \bo{0}) \, ,    
$$        
and we can write   
\begin{eqnarray}  \label{eq:appdifker}   
& & \mathcal{D}_k(F_{T_u})- \mathcal{D}_k(F_{u}) 
=  \normh{\esp\left[k\left(\cdot, \BG^{fo}_u \right)\right]}^2 - \normh{\esp\left[k\left(\cdot, \BG^{tot}_u \right)\right]}^2  \nonumber \\     
& \geq &    2\left<\esp\left[k\left(\cdot, \BG^{tot}_u\right)\right], \; 
 \esp\left[k\left(\cdot, \BG^{fo}_u \right)\right]- \esp\left[k\left(\cdot, \BG^{tot}_u \right)\right] \right>_{\mathcal{H}}   \nonumber \\  
	& =&   2  \esp\left[k\left(\BG^{fo}_u, \,  \BG^{tot '}_u \right)\right]
- 2  \esp\left[k\left(\BG^{tot}_u, \,  \BG^{tot '}_u\right)\right]   \nonumber  \, .     
\end{eqnarray}        
Using (\ref{eq:linktotfo}), Point (iii) holds by applying the Jensen inequality to $-k$, which is convex.  
  
\section{Proof of Theorem \ref{theo:sckerg}} \label{app:theo:sckerg}
Without loss of generality, we suppose that the outputs $\bo{Y} := \MG(\bo{\X}_{w}, \bo{Z}_{\sim w})$ is centered, that is, $\esp\left[\bo{Y} \right]=\bo{0}$. Recall that AFs are also centered. Using $w \subseteq u$, we can write $u=w \cup w_0$ with  $w_0 \subseteq u$ and $w \cap w_0 =\emptyset$. Thus, $\bo{Y} := \MG(\bo{\X}_{w}, \bo{Z}_{w_0}, \bo{Z}_{\sim u})$. \\  
First, as $
\MG^{fo}_w (\bo{\X}_{w}) =\esp_{\bo{Z}_{w_0} \bo{Z}_{\sim u}}\left[ \MG(\bo{\X}_{w}, \bo{Z}_{w_0}, \bo{Z}_{\sim u})\right]$,   
and it is known that  (see \cite{lamboni21ar}; Lemma 3)  
$  
\MG^{fo}_u (\bo{\X}_{u}) \stackrel{d}{=}  \MG^{fo}_u (\bo{\X}_{w}, \bo{Z}_{w_0})  =
\esp_{\bo{Z}_{\sim u}}\left[ \MG(\bo{\X}_{w}, \bo{Z}_{w_0}, \bo{Z}_{\sim u})\right]
$, we can see that        
$$  
\MG^{fo}_w (\bo{\X}_{w}) \stackrel{d}{=} \esp_{\bo{Z}_{w_0}}\left[\MG^{fo}_u (\bo{\X}_{w}, \bo{Z}_{w_0}) \right] \, .     
$$         
Second, for the convex kernel $k$, the Jensen inequality allows for writting $\mathcal{D}_k(F_{w})$ as      
$$
\mathcal{D}_k(F_{w}) =   \esp\left[k\left(\MG^{fo}_w (\bo{\X}_{w}), \,  \MG^{fo}_w (\bo{\X}_{w}') \right)\right] -k(\bo{0}, \bo{0})  \, .
$$  
Thus, the first result holds by applying the Jensen inequality, that is, 
$$   
 \esp\left[k\left(\MG^{fo}_w (\bo{\X}_{w}), \,  \MG^{fo}_w (\bo{\X}_{w}') \right)\right]   
\leq   
 \esp\left[k\left(\MG^{fo}_u (\bo{\X}_{u}), \,  \MG^{fo}_u (\bo{\X}_{u}') \right)\right] \, . 
$$  
  
For the second result, it comes out from the above equivalent in distribution that 
$$
\MG^{tot}_w (\bo{\X}_{w}, \bo{Z}_{w_0}, \bo{Z}_{\sim u}) 
= \MG (\bo{\X}_{w}, \bo{Z}_{w_0}, \bo{Z}_{\sim u}) - 
\esp_{\bo{\X}_{w}} \left[\MG (\bo{\X}_{w}, \bo{Z}_{w_0}, \bo{Z}_{\sim u}) \right] \, , 
$$  
$$
\MG^{tot}_u (\bo{\X}_{w}, \bo{Z}_{w_0}, \bo{Z}_{\sim u}) 
= \MG (\bo{\X}_{w}, \bo{Z}_{w_0}, \bo{Z}_{\sim u}) - 
\esp_{\bo{\X}_{w}, \bo{Z}_{w_0}} \left[\MG (\bo{\X}_{w}, \bo{Z}_{w_0}, \bo{Z}_{\sim u}) \right] \, , 
$$
and we want to show that 
$$     
\esp\left[k\left(\MG^{tot}_w (\bo{\X}_{w}, \bo{Z}_{\sim w}), \,  \MG^{tot}_w (\bo{\X}_{w}', \bo{Z}'_{\sim w}) \right)\right]    
\leq  
 \esp\left[k\left(\MG^{tot}_u (\bo{\X}_{u}, \bo{Z}_{\sim u}), \,  \MG^{tot}_u (\bo{\X}_{u}', \bo{Z}'_{\sim u}) \right)\right] \, . 
$$    
To that end, let  $\bo{V}' :=(\bo{\X}_{w}', \, \bo{Z}_{w_0}', \, \bo{Z}_{\sim u}')$ be an i.i.d. copy of $\bo{V} := (\bo{\X}_{w}, \bo{Z}_{w_0}, \bo{Z}_{\sim u})$; 
and consider the function
$   
h(\bo{V}, \bo{V}') := \MG (\bo{\X}_{w}, \bo{Z}_{w_0}, \bo{Z}_{\sim u}) - 
\MG (\bo{\X}_{w}', \bo{Z}_{w_0}', \bo{Z}_{\sim u}') 
$. 
  Since the three components of $\bo{V}$ (resp.  $\bo{V}'$) are independent, we can write 
$$
\esp\left[h(\bo{V}, \bo{V}') \, | \bo{\X}_{w},\, \delta_{\bo{0}}(\bo{Z}_{w_0}' -\bo{Z}_{w_0}),\, \delta_{\bo{0}}(\bo{Z}_{\sim u}'-\bo{Z}_{\sim u})  \right] =  \MG^{tot}_w (\bo{\X}_{w}, \bo{Z}_{w_0}, \bo{Z}_{\sim u})  \, ,
$$ 
$$
\esp\left[h(\bo{V}, \bo{V}') \, | \bo{\X}_{w}, \bo{Z}_{w_0}, \delta_{\bo{0}}(\bo{Z}_{\sim u}'-\bo{Z}_{\sim u})   \right]   
=  \MG^{tot}_u (\bo{\X}_{w}, \bo{Z}_{w_0}, \bo{Z}_{\sim u})  \, .
$$   
Moreover,  the properties of conditional expectation allow for writing 
\begin{eqnarray}
&& \esp\left[ \MG^{tot}_u (\bo{\X}_{w}, \bo{Z}_{w_0}, \bo{Z}_{\sim u})\, |\bo{\X}_{w},\, \delta_{\bo{0}}(\bo{Z}_{w_0}' -\bo{Z}_{w_0}),\, \delta_{\bo{0}}(\bo{Z}_{\sim u}'-\bo{Z}_{\sim u})  \right] \nonumber \\ 
&=& \esp\left[  
\esp\left[h(\bo{V}, \bo{V}') \, |\bo{\X}_{w}, \bo{Z}_{w_0}, \delta_{\bo{0}}(\bo{Z}_{\sim u}'-\bo{Z}_{\sim u})   \right]   \, | \bo{\X}_{w},\, \delta_{\bo{0}}(\bo{Z}_{w_0}' -\bo{Z}_{w_0}),\, \delta_{\bo{0}}(\bo{Z}_{\sim u}'-\bo{Z}_{\sim u})  \right]   \nonumber \\ 
 &=&  \esp\left[h(\bo{V}, \bo{V}') \, | \bo{\X}_{w},\, \delta_{\bo{0}}(\bo{Z}_{w_0}' -\bo{Z}_{w_0}),\, \delta_{\bo{0}}(\bo{Z}_{\sim u}'-\bo{Z}_{\sim u})  \right] 
=  \MG^{tot}_w (\bo{\X}_{w}, \bo{Z}_{w_0}, \bo{Z}_{\sim u})   \, , \nonumber 
\end{eqnarray}   
because the space of projection and the filtration associated with $(\bo{\X}_{w}, \bo{Z}_{w_0}, \, \delta_{\bo{0}}(\bo{Z}_{\sim u}'-\bo{Z}_{\sim u}) )$ contain those of $(\bo{\X}_{w},\, \delta_{\bo{0}}(\bo{Z}_{w_0}' -\bo{Z}_{w_0}),\, \delta_{\bo{0}}(\bo{Z}_{\sim u}'-\bo{Z}_{\sim u}))$.  The second result holds by applying the conditional Jensen inequality, as $k$ is convex. \\ 
Finally, the results for a concave kernel $k$ can be deduced from the above results. Indeed, we can see that $-k$ is convex and $\mathcal{D}_k(F_{w})$ becomes  
$$
\mathcal{D}_k(F_{w}) =   \esp\left[-k\left(\MG^{fo}_w (\bo{\X}_{w}), \,  \MG^{fo}_w (\bo{\X}_{w}') \right)\right] + k(\bo{0}, \bo{0})  \, .
$$

\section{Proof of Corollary \ref{coro:proimnor}} \label{app:coro:proimnor}
It is sufficient to show the results for $q=1$.\\ 
For Point (i), according to Theorem \ref{theo:sckerg}, we can write 
$$
0\leq \mathcal{D}_k(F_{w})  \leq   \mathcal{D}_k(F_{T_w})  \leq   \mathcal{D}_k(F_{T_u}), \quad \forall\, w \subseteq u \subseteq \{1, \ldots, d\} \, . 
$$ 
Thus, we have $0\leq S_k (F) \leq 1$ because $F_\bu =F_{T_{\{1, \ldots, d\}}}$. \\  
Point (ii) is obvious because $k\in \mathcal{K}_E$, the set of kernels that guarantee the independence criterion. \\  
The if part of Point (iii) is obvious. For the only if part, the equality $\mathcal{D}_k(F_{T_u})  =   \mathcal{D}_k(F_{\bu})$ implies that  
$$    
\int k(\bo{y}, \bo{y}') \, d(F_{T_u} \otimes F_{T_u} -F_{\bu} \otimes F_{\bu})(\bo{y}, \bo{y}') = 0\, ,       
$$    
which also implies that $F_{T_u} = F_{\bu}$ for the second kind of kernels of $\mathcal{K}_E$. \\
Point (iv) holds for IMKs by definition.

\section{Proof of Theorem \ref{theo:teststat}} \label{app:theo:teststat}
Firstly, we have $\widehat{\mu}(\bo{Z}_i)- \esp_{\bo{\X}_u}\left[\MG(\bo{\X}_u , \bo{Z}_i) \right] \to 0$ when $m_1 \to  \infty$. \\
Knowing that 
$ 
\BG_{i,u}^{tot} = \MG(\bo{\X}_{i,u} , \bo{Z}_i) - \esp_{\bo{\X}_u} \left[\MG(\bo{\X}_{u} , \bo{Z}_i) \right] 
$ 
and 
$
\BG_{i,u}^{tot\, '} = \MG(\bo{\X}_{i,u}' , \bo{Z}_i') - \esp_{\bo{\X}_u'} \left[\MG(\bo{\X}_{u}' , \bo{Z}_i') \right]   
$,
the Taylor expansion of $k$ about $\left(\BG_{i,u}^{tot} ,\, \BG_{i,u}^{tot\, '}\right)$ yields 
\begin{eqnarray}   
& & k\left(\MG(\bo{\X}_{i,u} , \bo{Z}_i) - \widehat{\mu}(\bo{Z}_i),\, \MG(\bo{\X}_{i,u}' , \bo{Z}_i') - \widehat{\mu}(\bo{Z}_i') \right) 
 =  k\left(\BG_{i,u}^{tot} ,\, \BG_{i,u}^{tot\, '}\right) \nonumber  \\ 
& &  + \nabla^\T k \left(\BG_{i,u}^{tot} ,\, \BG_{i,u}^{tot\, '}\right) \left[\begin{array}{c}    \esp_{\bo{\X}_u}\left[\MG(\bo{\X}_u , \bo{Z}_i) \right] - \widehat{\mu}(\bo{Z}_i) \\
 \esp_{\bo{\X}_u'}\left[\MG(\bo{\X}_u' , \bo{Z}_i') \right]   - \widehat{\mu}(\bo{Z}_i') \end{array}\right] +  R_{m_1} \nonumber  \, ,
\end{eqnarray}  
where $R_{m_1} \xrightarrow{P} 0$ when $m_1 \to~\infty$. Therefore, we can write   
\begin{eqnarray}
 & & \widehat{\mu_k^{tot}}  
 =  \frac{1}{m} \sum_{i=1}^m k\left(\BG_{i,u}^{tot} ,\, \BG_{i,u}^{tot\, '}\right)  \nonumber  \\ 
& &  + \frac{1}{m} \sum_{i=1}^m  \nabla^\T k \left(\BG_{i,u}^{tot} ,\, \BG_{i,u}^{tot\, '}\right) \left[\begin{array}{c}  \esp_{\bo{\X}_u}\left[\MG(\bo{\X}_u , \bo{Z}_i) \right] - \widehat{\mu}(\bo{Z}_i) \\  \esp_{\bo{\X}_u'}\left[\MG(\bo{\X}_u' , \bo{Z}_i') \right]   - \widehat{\mu}(\bo{Z}_i') 
\end{array}\right]  + R_{m,m_1} \nonumber  \, ,
\end{eqnarray}   
where $R_{m,m_1} \xrightarrow{P} 0$ when $m_1 \to~\infty$. Since the second term of the above equation converge in probability toward $0$, the LLN ensures that $\widehat{\mu_k^{tot}}$ is a consistent estimator of $\mu_k^{tot}$. thus, the first result of Point (i) holds.\\    
Secondly, we obtain the second result of Point (i) by applying the central limit theorem (CLT) to the first term of the above equation, as the second term converge in probability toward $0$.  \\         
The proof of Point (ii) is similar to the proof of Point (i). Indeed,  using the Taylor expansion of $k^2$, we obtain the consistency of the second-order moment of $k$. The Slutsky theorem ensures the consistency of the cross components and $\left( \widehat{\mu_k^{tot}}\right)^2$.  \\   
Point (iii) is then obvious using Point (ii). \\  
The proofs of Point (iv) is similar to those of  Point (i).

\section{Proof of Corollary \ref{coro:estksi}} \label{app:coro:estksi}
First, the results about the consistency of the estimators are obtained by using Theorem \ref{theo:teststat} and the Slutsky theorem. \\  
The numerators of Equations (\ref{eq:estksifo})-(\ref{eq:estksitot}) are asymptotically distributed as Gaussian variable according to Theorem \ref{theo:teststat}. To obtain the asymptotic distributions of the sensitivity indices, we first applied  the Slutsky theorem, and second, we use the fact that $\sqrt{m}\left(\widehat{S_{T_u}^k} - \frac{\esp\left[k\left(\MG^{tot}_u,\, \MG^{tot\, '}_u \right) \right] -k(\bo{0}, \bo{0})}{\frac{1}{M}\sum_{i=1}^M k\left(\MG(\bo{\X}_{i,u} , \bo{Z}_i) - \widehat{\mu},\, \MG(\bo{\X}_{i,u}' , \bo{Z}_i') - \widehat{\mu} \right) -k(\bo{0}, \bo{0})} \right)$ and $\sqrt{m}\left(\widehat{S_{T_u}^k} - S_{T_u}^k \right)$ are asymptotically  equivalent in probability under the technical condition $m / M \to 0$ (see \cite{lamboni18} for more details).        
    
\section{Proof of Lemma \ref{lem:conker}} \label{app:lem:conker}
For Point (i), the convexity of $\psi$ implies the existence of $\partial \psi$ such that  
$$
-\alpha \psi(\bo{y}, \bo{y}') +\alpha \psi(\bo{b}, \bo{y}') \leq \left< -\alpha \partial \psi(\bo{b}, \bo{y}'), \bo{y}-\bo{b} \right> \, ,
$$      
which also implies  (thanks to the Taylor expansion) that
\begin{equation} \label{eq:psicon}
e^{-\alpha\psi(\bo{y}, \bo{y}') + \alpha \psi(\bo{b}, \bo{y}')} \leq 
 e^{\left< -\alpha \partial \psi(\bo{b}, \bo{y}'), \bo{y}-\bo{b} \right>}
\approx 
1 + \left< -\alpha \partial \psi(\bo{b}, \bo{y}'), \bo{y}-\bo{b} \right> \,  
\end{equation}
under the condition (thanks to Cauchy-Schwartz) 
$$ 
\left| \left< -\alpha \partial \psi(\bo{b}, \bo{y}'), \bo{y}-\bo{b} \right> \right| \leq 
\alpha \norme{\partial \psi(\bo{b}, \bo{y}')} \, \norme{\bo{y}-\bo{b}} \leq \epsilon
\qquad  \forall\, \bo{y}', \bo{b} \in \R^\NN \, . 
$$          
 Equivalently, we can write  
$
 \alpha 
 \leq   \frac{\epsilon}{\norme{\bo{y}-\bo{b}} \norme{\partial \psi(\bo{b}, \bo{z})}}; 
\quad \forall\,  \, \bo{y}, \bo{z}, \bo{b} \in \mathcal{\X}    
$.  
Equation (\ref{eq:psicon}) implies that $k$ is concave under the above condition. Indeed, we have
\begin{eqnarray} 
e^{-\alpha\psi(\bo{y}, \bo{y}') + \alpha \psi(\bo{b}, \bo{y}')} -1 &\leq & \left< -\alpha \partial \psi(\bo{b}, \bo{y}'), \bo{y}-\bo{b} \right> \nonumber  \\
e^{-\alpha\psi(\bo{y}, \bo{y}')} - e^{-\alpha \psi(\bo{b}, \bo{y}')}  &\leq &  e^{-\alpha \psi(\bo{b}, \bo{y}')} \left< -\alpha \partial \psi(\bo{b}, \bo{y}'), \bo{y}-\bo{b} \right>
\nonumber \\  
- k (\bo{y}, \bo{y}') + k(\bo{b}, \bo{y}') &\geq & \left<\alpha \partial \psi(\bo{b}, \bo{y}') k(\bo{b}, \bo{y}'), \bo{y}-\bo{b} \right> 
= \left< \partial k(\bo{b}, \bo{y}'), \bo{y}-\bo{b} \right> 
\nonumber \, ,      
\end{eqnarray}   
with $\partial k(\bo{b}, \bo{y}') := \alpha \partial \psi(\bo{b}, \bo{y}') k(\bo{b}, \bo{y}')$ the subgradient of $-k$. Thus, $-k$ is convex because $k$ is continuous (\cite{boyd04}). \\ 
For Point (ii), the gradient and the hessian of  $k (\bo{y}, \bo{y}')=e^{- \alpha \psi(\bo{y}, \bo{y}')}$ w.r.t. $\bo{y}$ are
$$ 
\nabla k (\bo{y}, \bo{y}')  := -\alpha \nabla \psi(\bo{y}, \bo{y}') k(\bo{y}, \bo{y}'),  \, ,  
$$
$$ 
H_k (\bo{y}, \bo{y}')  := \left[ -\alpha H_\psi(\bo{y}, \bo{y}') +
\alpha^2 \nabla \psi(\bo{y}, \bo{y}') \nabla^\T \psi(\bo{y}, \bo{y}')  \right] k(\bo{y}, \bo{y}')  \, . 
$$  
Therefore, if we use $E := -H_\psi(\bo{y}, \bo{y}') + \alpha \nabla \psi(\bo{y}, \bo{y}') \nabla^\T \psi(\bo{y}, \bo{y}')$, then $k$ is concave when 
 $E$ is negative definite. Thus, for all $\bo{b} \in \mathcal{\X}$, we can write 
\begin{eqnarray}   
- \bo{b}^\T  H_\psi(\bo{y}, \bo{y}') \bo{b} +
\alpha  \bo{b}^\T \nabla \psi(\bo{y}, \bo{y}') \nabla^\T \psi(\bo{y}, \bo{y}') \bo{b} &\leq 0 & \nonumber \\
\alpha \left( \bo{b}^\T \nabla \psi(\bo{y}, \bo{y}')\right)^2 
& \leq & \bo{b}^\T  H_\psi(\bo{y}, \bo{y}') \bo{b} \nonumber  \\
\alpha &\leq & 
\frac{\bo{b}^\T  H_\psi(\bo{y}, \bo{y}') \bo{b} }{\left( \bo{b}^\T \nabla \psi(\bo{y}, \bo{y}')\right)^2} \nonumber  \, .   
\end{eqnarray}

\section{Proof of Corollary \ref{coro:conkerunb}} \label{app:coro:conkerunb}
Namely, we use $u_1(\bo{y}, \bo{y}', \bo{y}'') :=  \frac{\epsilon}{\norme{\bo{y}-\bo{y}'} \norme{\partial \psi(\bo{y}', \bo{y}'')}}$ and $u_2(\bo{y}, \bo{y}', \bo{y}'') := \frac{\bo{y}^{'' \,\T}  H_\psi(\bo{y}, \bo{y}') \bo{y}''}{\left(\nabla \psi(\bo{y}, \bo{y}')^\T \, \bo{y}''\right)^2}$ for the upper bounds of $\alpha$ (see proof of  Lemma \ref{lem:conker}). For the sequel of simplicity, we use $u(\bo{y}, \bo{y}', \bo{y}'')$ with $\bo{y}, \bo{y}', \bo{y}'' \in \mathcal{\X}$ for either $u_1(\bo{y}, \bo{y}', \bo{y}'')$ or $u_2(\bo{y}, \bo{y}', \bo{y}'')$. \\    
As $u(\bo{Y}, \bo{Y}', \bo{Y}'')$ is random variable, we have (Markov's inequality)  
$$ 
\proba\left(u(\bo{Y}, \bo{Y}', \bo{Y}'') <  \alpha \right) = 
\proba\left(\frac{1}{u(\bo{Y}, \bo{Y}', \bo{Y}'')} > \frac{1}{\alpha} \right) 
\leq \alpha \esp\left[\frac{1}{u(\bo{Y}, \bo{Y}', \bo{Y}'')} \right] \leq \tau \, ,  
$$  
which implies that $\alpha  \leq \frac{\tau}{\esp\left[\frac{1}{u(\bo{Y}, \bo{Y}', \bo{Y}'')} \right]}$.   \\  
Moreover, using Markov's inequality we can write
$$      
1- \tau \leq  
\proba\left(u(\bo{Y}, \bo{Y}', \bo{Y}'') \geq  \alpha \right) 
\leq  \frac{1}{\alpha} \esp\left[u(\bo{Y}, \bo{Y}', \bo{Y}'') \right] \, , 
$$  
which implies that $\alpha  \leq \frac{\esp\left[u(\bo{Y}, \bo{Y}', \bo{Y}'') \right]}{1- \tau }$.

 \end{appendices}

	 
                 


\end{document}